\documentclass[graybox]{svmult}

\usepackage[numbers, square]{natbib}
\usepackage{type1cm}        
\usepackage{makeidx}         
\usepackage{multicol}        
\usepackage[bottom]{footmisc}

\usepackage{newtxtext}       %
\usepackage[varvw]{newtxmath}       

\usepackage{graphicx}
\usepackage{epstopdf}
\usepackage{algorithm}
\usepackage{hyperref}
\newcommand{\bq}{\begin{eqnarray*}}
\newcommand{\eq}{\end{eqnarray*}}
\newcommand{\bqn}{\begin{eqnarray}}
\newcommand{\eqn}{\end{eqnarray}}

\makeindex             

\begin{document}

\title*{Topological  State-Space Estimation of Functional Human Brain Networks}
\titlerunning{Topological  State-Space Estimation of Functional Human Brain Networks} 
\author{Moo K. Chung$^{*1}$, Shih-Gu Huang$^2$, Ian C. Carroll$^3$,Vince D. Calhoun$^4$, H. Hill Goldsmith$^5$}
\authorrunning{Chung et al.} 
\institute{$^{*1}$Moo K. Chung \at Department of Biostatistics and Medical Informatics, University of Wisconsin, Madison, WI, USA \email{mkchung@wisc.edu}
\and $^2$Shih-Gu Huang \at PPG/ECG Signal, Taiwan \email{shihgu@gmail.com}
\and $^3$Ian C. Carroll \at Department of Child and Adolescent Psychiatry, New York University Grossman School of Medicine, USA 
\email{ian.carroll@nyulangone.org}
 \and $^4$Vince D. Calhoun \at Tri-Institutional Center for Translational Research in Neuroimaging and Data Science (TReNDS), Georgia State, Georgia Tech, Emory Georgia State University, Georgia, USA
 \email{vcalhoun@gsu.edu}
 \and $^5$H. Hill Goldsmith \at Department of Psychology \& Waisman Center, University of Wisconsin-Madison, USA
 \email{hill.goldsmith@wisc.edu}
}

%
%
\maketitle

\abstract{
We introduce an innovative, data-driven topological data analysis (TDA) technique for estimating the state spaces of dynamically changing functional human brain networks at rest. Our method utilizes the Wasserstein distance to measure topological differences, enabling the clustering of brain networks into distinct topological states. This technique outperforms the commonly used k-means clustering in identifying brain network state spaces by effectively incorporating the temporal dynamics of the data without the need for explicit model specification. We further investigate the genetic underpinnings of these topological features using a twin study design,  examining the heritability of such state changes. Our findings suggest that the topology of brain networks, particularly in their dynamic state changes, may hold significant hidden genetic information. MATLAB code for the method is available at \url{https://github.com/laplcebeltrami/PH-STAT}.
}

\section{Author Summary}

The paper introduces a new data-driven topological data analysis (TDA) method for studying dynamically changing human functional brain networks obtained from the resting-state functional magnetic resonance imaging (rs-fMRI). Leveraging persistent homology, a  multiscale topological approach, we present a framework that incorporates the temporal dimension of brain network data. This allows for a more robust estimation of the topological features of dynamic brain networks.

The method employs the Wasserstein distance to measure the topological differences between networks and demonstrates greater efficiency and performance than the commonly used $k$-means clustering in defining the state spaces of dynamic brain networks. 
Our method maintains robust performance across different scales and is especially suited for dynamic brain networks.

In addition to the methodological advancement, the paper applies the proposed technique to analyze the heritability of overall brain network topology using a twin study design. The study investigates whether the dynamic pattern of brain networks is a genetically influenced trait, an area previously underexplored. By examining the state change patterns in twin brain networks, we make significant strides in understanding the genetic factors underlying dynamic brain network features. Furthermore, the paper makes its method accessible by providing MATLAB codes, contributing to reproducibility and broader application.

\section{Introduction}
\label{sec:Introduction}

In standard graph theory-based network analysis, network features such as node degrees and clustering coefficients are obtained from adjacency matrices after thresholding weighted edges \citep{bassett.2017,sporns.2003,vanwijk.2010,chung.2017.BC}. The final statistical analysis results can vary depending on the choice of threshold or parameter \citep{chung.2013.MICCAI,lee.2012.tmi}. This variability underscores the need for a multiscale network analysis framework that provides consistent results and interpretation, regardless of the choice of parameter. Persistent homology, a branch of algebraic topology, presents a novel solution to this challenge of multiscale analysis \citep{edelsbrunner.2010}. Unlike traditional graph theory approaches that analyze networks at a single fixed scale, persistent homology examines networks across multiple scales. It identifies topological features that remain persistent and are robust against different scales and noise perturbations 
\citep{petri.2014, sizemore.2018, sizemore.2019,vaccarino.2022}. 

Recent studies have  illustrated the versatility of persistent homology in analyzing complex networks, including brain networks. \citet{sizemore.2019,xing.2022} highlighted the application of persistent homology in evaluating temporal changes in topological network features. 
\citet{aktas.2019} used persistent homology to detect and track the evolution of networks’ clique. \citet{billings.2021} discussed the use of simplicial complexes encoded by persistent homology for brain networks. \citet{sizemore.2018} applied persistent homology to investigate the spatial distributions of cliques and cycles in brain networks. \citet{khalid.2014,caputi.2021} showed how persistent homology could be used in the analysis of functional brain connectivity using EEG. \citet{chung.2023.NI} utilized persistent homology to analyze brain networks for studying abnormal white matter in maltreated children. These studies collectively emphasize the potential of persistent homology in providing a robust framework for multiscale network analysis. This approach's ability to capture topological features across different scales and under varying conditions makes it particularly suitable for studying the complex brain networks.

Persistent homological network approaches have shown to be more robust and outperform many existing graph theory measures and methods. 
In \citet{lee.2011.MICCAI, lee.2012.tmi}, persistent homology was shown to outperform eight existing graph theory features, such as clustering coefficient, small-worldness, and modularity. \citet{kuang.2019} showed persistent homology-based measures can provide  more significant group difference and better classification performance compared to  standard graph-based measures that characterize small-world organization and modular structure.  In \citet{chung.2017.CNI,chung.2019.ISBI}, persistent homology was shown to outperform various matrix norm-based network distances. In \citet{wang.2018.annals}, persistent homology was shown to outperform the power spectral density and local variance methods. In \citet{wang.2017.CNI}, persistent homology was shown to outperform topographic power maps. In \citep{yoo.2017}, center persistency was shown to outperform the network-based statistic and element-wise multiple corrections. In 
\citet{chung.2023.wasserstein}, persistent homology based clustering is shown to outperform $k$-means clustering and hierarchical clustering. 
However, the method has been mainly used on {\em static} networks or a static summary of time-varying networks. The dynamic pattern of persistent homology for time-varying brain networks was rarely investigated, with a few exceptions \citep{yoo.2016, santos.2019,song.2020.ISBI, giusti.2016, sizemore.2018,chung.2023.wasserstein}.

While Euclidean loss remains the dominant cost function in deep learning, topological losses based on persistent homology are emerging as superior in tasks requiring topological understanding \citep{chen.2019, hu.2019, gupta.2022, lin.2023}. These topological losses incorporate penalties based on the topological features of the data, distinguishing them from the Euclidean loss, which primarily focuses on differences at the node or edge level. By encoding the intrinsic topological structure of the network, topological losses facilitate the creation of more informative feature maps, potentially enhancing overall model performance \citep{hofer.2019}. \citet{gupta.2022} demonstrated that image segmentation based on topological loss outperforms other deep learning architectures for similar tasks. \citet{lin.2023} introduced a new architecture that excels in segmenting curvilinear structures by learning topological similarities over existing methods.

In this paper, we propose to develop a novel {\em dynamic persistent homology} framework for time varying network data. Our coherent scalable framework for the computation is based on the Wasserstein distance between persistent diagrams, which provides the topological profile of data into 2D scatter plots. We directly establish the relationship between the Wasserstein distance and edge weights in networks  making the method far more accessible and adaptable. 
We achieve $\mathcal{O}(n \log n)$ run time in most graph manipulation tasks such as matching and averaging. Such scalable computation enables us to perform a computationally demanding task such as topological clustering with ease. The method is applied in the determination of the  state space of dynamically changing functional brain networks obtained from the resting-state functional magnetic resonance imaging (rs-fMRI).  We will show that  the proposed  method based on the Wasserstein  distance can capture the topological patterns that are consistently observed  across different time points.

The Wasserstein distance or Kantorovich–Rubinstein metric, as  originally defined between probability distributions, can be used to measure topological differences \citep{vallender.1974,canas.2012,berwald.2018}. Due to the connection to the optimal mass transport, which enjoys  various optimal properties, the Wasserstein distance has been applied to various imaging applications. Nonetheless, its application in network data analysis remains relatively limited \citep{ma.2023,chung.2023.NI}.
\citet{mi.2018} used the Wasserstein distance in resampling brain surface meshes. \citet{shi.2016,su.2015} used the Wasserstein distance in classifying brain cortical surface shapes. \citet{hartmann.2018} used  the Wasserstein distance in building generative adversarial networks. \citet{sabbagh.2019} used the Wasserstein distance for manifold regression problems in the space of positive definite matrices for the source localization problem in EEG. \citet{xu.2021} used the Wasserstein distance in predicting Alzheimer's disease progression in magnetoencephalography (MEG) brain networks. \citet{fu.2023} enhanced images by regularizing with the Wasserstein distance. However,  the Wasserstein distance in these applications is all geometric in nature.  

We applied the method to dynamically changing twin brain networks obtained from the resting-state functional magnetic resonance imaging (rs-fMRI). We investigated if the state change pattern in time varying brain networks is genetically heritable for the first time. This is not yet reported in existing literature. Monozygotic (MZ) twins share 100\% of genes while  dizygotic (DZ) twins share 50\% of genes \citep{falconer.1995}. MZ-twins are more similar or concordant than DZ-twins for cognitive aging and  dysfunction \citep{reynolds.2015}. The difference between MZ- and DZ-twins directly quantifies the extent to which imaging phenotypes, behaviors and cognitions are influenced by genetic factors \citep{zhan.2022}. If MZ-twins show more similarity on a given trait compared to DZ-twins, this provides a piece of evidence that genes significantly influence that trait. 
Even twin studies on normal subjects are useful for understanding the extent to which psychological and medical disorders, as well as behaviors and traits, are influenced by genetic factors. This information can be used to develop better ways to prevent and treat disorders and maladaptive behaviors.  Some of the most effective treatments for medical disorders have been identified as a result of twin studies \citep{sahu.2016}. 

Even though there are numerous twin imaging studies, almost all previous studies used {\em static} univariate imaging phenotypes such as cortical thickness \citep{mckay.2014}, fractional anisotropy \citep{chiang.2011}, functional activation \citep{blokland.2011,glahn.2010,smit.2008} in  determining heritability in brain networks.  There have been a limited number of  studies investigating the heritability of  the {\it dynamics} of brain networks \citep{blokland.2011,vidaurre.2017}. {\it It is  not even clear the dynamic pattern itself is a heritable trait.} We propose to tackle this challenge. Measures of network dynamics are worth investigating as potential phenotypes that indicate the genetic risk for  neuropsychiatric disorders \citep{bullmore.2009}. Determining the extent of  heritability of dynamic pattern is the first necessary prerequisite for identifying dynamic network phenotypes. 

One of the  earliest papers on functional brain activation in twins is based on the resting-state EEG \citep{lykken.1982}, where they observed high twin correlation in MZ-twins on EEG spectra. \citet{glahn.2010} reported a heritability of 0.42  for default-mode network (DMN) in an extended 
pedigree study without twins. \citet{xu.2017} reported a heritability of 0.54 for DMN on using 24 pairs of MZ and 22 pairs of DZ. \citet{korgaonkar.2014} studied 79 MZ twins and 46 DZ twin pairs, reporting heritability in only one specific connection: They found statistically significant heritability of 0.41 for the connection between the precuneus and the right inferior parietal/temporal cortex, using  a structural equation model. We report far stronger results with much higher heritability in a larger twin study.

\section{Methods}

\begin{figure}[t]
\centering
\includegraphics[width=1\linewidth]{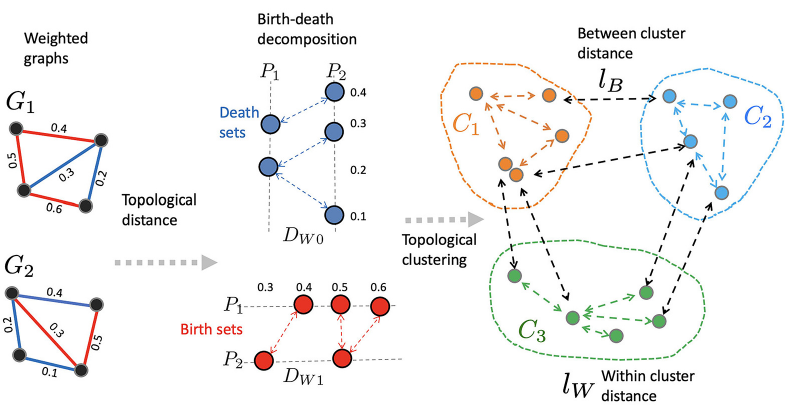}
\caption{Proposed topological clustering pipeline used in estimating the state space. Given two weighted graphs \( G_1, G_2 \), we first perform the birth-death decomposition and partition the edges into sorted birth and death sets (section \ref{sec:BDD}). 
The 0D topological distance $D_{W0}$ between birth values quantifies discrepancies in connected components (section \ref{sec:dist}). 
The 1D topological distance $D_{W1}$ between death values quantifies discrepancies in cycles (section \ref{sec:dist}). The combined distance $\mathcal{D}=D_{W0}^2 + D_{W1}^2$ is used in computing the within-cluster distance $l_W$ between graphs. Topological clustering is performed by minimizing $l_W$ over all possible cluster labels $C_1, \cdots, C_k$ (section \ref{sec:clustering}).}
\label{fig:schematic1}
\end{figure}

\subsection{Graphs as simplicial complices}
The proposed method for estimating topological state space is based on the topological clustering on a collection of graphs (Figure \ref{fig:schematic1}). The initial step involves a birth-death decomposition of a weighted graph, leading to the generation of sorted birth and death sets (section \ref{sec:BDD}). The second step entails calculating the topological distance: between birth sets to obtain the 0D topological distance \(d_0\), and between death sets to obtain the 1D topological distance \(d_1\) (section \ref{sec:dist}). The third step involves computing the within-cluster distance \(l_W\) among the collection of graphs (section \ref{sec:clustering}). Subsequently, we demonstrate the equivalence of topological clustering with \(k\)-means clustering in a high-dimensional convex set, employing \(k\)-means clustering routines for optimization.
To increase the reproducibility, MATLAB codes for performing the methods are provided in \url{https://github.com/laplcebeltrami/PH-STAT}.

A high dimensional object such as a brain network can be modeled as weighted graph $\mathcal{X} = (V, w)$ consisting of node set $V$ indexed as $V=\{1, 2, \cdots,$ $p\}$ and edge weights $w=(w_{ij})$ between nodes $i$ and $j$.  If we order the edge weights in the increasing order, we have the sorted edge weights:
\bqn \min_{j,k} w_{jk} = w_{(1)} < w_{(2)} < \cdots < w_{(q)} = \max_{j,k} w_{jk}, \eqn
where $q \leq (p^2-p)/2$.  The subscript $_{( \;)}$ denotes the order statistic. In terms of sorted edge weight set
$W=\{ w_{(1)}, \cdots, w_{(q)} \},$ we may also write the graph as $\mathcal{X} = (V, W)$. If we connect nodes following some criterion on the edge weights, they will form a simplicial complex which will follow the topological structure of the underlying weighted graph \citep{edelsbrunner.2010,zomorodian.2009}.  Note that the $k$-simplex  is the convex hull of $k+1$ points in $V$.   A  simplicial complex is a finite collection of simplices such as points (0-simplices), lines (1-simplices), triangles (2-simplices) and higher dimensional counter parts.

The {\em Rips complex} $\mathcal{X}_{\epsilon}$ is a simplicial complex, whose $k$-simplices are formed by $(k+1)$ nodes which are pairwise within distance $\epsilon$ \citep{ghrist.2008}. While a graph has at most 1-simplices, the Rips complex has at most $(p-1)$-simplices. The Rips complex induces a hierarchical nesting structure called  the Rips  filtration
$$\mathcal{X}_{\epsilon_0} \subset \mathcal{X}_{\epsilon_1}\subset \mathcal{X}_{\epsilon_2}\subset \cdots $$ 
for $0=\epsilon_{0} < \epsilon_{1} < \epsilon_{2} < \cdots$, where the sequence of $\epsilon$-values are called the filtration values. The filtration is quantified through a topological basis called {\em $k$-cycles}. 0-cycles are the connected components, 1-cycles are 1D closed paths or loops while 2-cycles are 3-simplices (tetrahedron) without interior. Any $k$-cycle can be represented as a linear combination of basis $k$-cycles. The Betti number $\beta_k$ counts the number of independent $k$-cycles. During the Rips filtration, the $i$-th $k$-cycle is born at filtration value $b_i$ and dies at $d_i$. The collection of all the paired filtration values 
$$P(\mathcal{X})=\{ (b_1, d_1), \cdots, (b_q, d_q) \}$$ displayed as 1D intervals is called the {\em barcode} and displayed as a 2D scatter plot is called the {\em persistent diagram}. Since $b_i < d_i$, the scatter plot in the persistent diagram  are displayed above the line $y=x$ line by taking births in the $x$-axis and deaths in the $y$-axis. 

For a dynamically changing brain network $\mathcal{X}(t) = (V, w(t))$, we assume the node set is fixed while edge weights are changing  over time $t$. If we build persistent homology at each fixed time, the resulting barcode is also time dependent:
 $$P(\mathcal{X}(t)) = \{ (b_1(t), d_1(t)), \cdots, (b_q(t), d_q(t)) \}.$$

\subsection{Graph filtrations}

As the number of nodes $p$ increases, the resulting Rips complex becomes increasingly dense. Additionally, as the filtration values rise, the number of edges connecting each pair of nodes also increases, leading to a more interconnected structure. At higher filtration values, Rips filtration becomes an ineffective representation of networks. To remedy this problem, graph filtration was introduced   \citep{lee.2011.MICCAI,lee.2012.tmi}. Given weighted graph $\mathcal{X}=(V, w)$ with edge weight $w = (w_{ij})$, the binary network $\mathcal{X}_{\epsilon} =(V, w_{\epsilon})$ is a graph consisting of the node set $V$ and the binary edge weights 
$w_{\epsilon} =(w_{\epsilon,ij})$ given by 
\bq w_{\epsilon,ij} =   \begin{cases}
1 &\; \mbox{  if } w_{ij} > \epsilon;\\
0 & \; \mbox{ otherwise}.
\end{cases}
\label{eq:case}
\eq
Note $w_{\epsilon}$ is the adjacency matrix of $\mathcal{X}_{\epsilon}$, which is a simplicial complex consisting of $0$-simplices (nodes) and $1$-simplices (edges)  \citep{ghrist.2008}. While the binary network $\mathcal{X}_{\epsilon}$ has at most 1-simplices, the Rips complex can have at most $(p-1)$-simplices. 
By choosing threshold values at sorted edge weights $w_{(1)}, w_{(2)}, \cdots, w_{(q)}$, we obtain the sequence of nested graphs \citep{chung.2013.MICCAI}:
\bqn \mathcal{X}_{w_{(1)}} \supset \mathcal{X}_{w_{(2)}} \supset \cdots \supset \mathcal{X}_{w_{(q)}}. \eqn
The sequence of such a nested multiscale graph  is called the {\em graph filtration} \citep{lee.2011.MICCAI,lee.2012.tmi}. 
Note that $\mathcal{X}_{w_{(1)} - \epsilon}$ is the complete weighted graph for any $\epsilon>0$. 
On the other hand, $\mathcal{X}_{w_{(q)}}$ is the node set $V$. By increasing the threshold value, we are thresholding at higher connectivity; thus more edges are removed.

\begin{figure}[t]
\begin{center}
\includegraphics[width=1\linewidth]{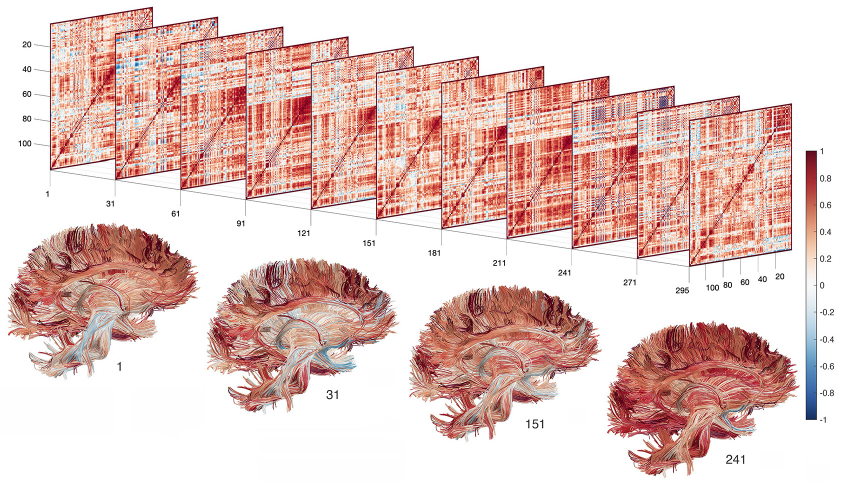}
\caption{Dynamically changing correlation matrices computed from rs-fMRI using the sliding window of size 60 for a subject. The constructed correlation matrices are superimposed on top of the white matter fibers in the MNI space and  color coded based on correlation values.}
\label{fig:dynamicTDA}
\end{center}
\end{figure}

For dynamically changing brain networks (Figure \ref{fig:dynamicTDA}), we  can similarly build time varying graph filtrations at each time point $\{\mathcal{X}_\epsilon(t): t \in \mathbb{R^+} \}$.

\subsection{Birth-death decomposition}
\label{sec:BDD}
Unlike the Rips complex, there are no higher dimensional topological features beyond the 0D and 1D topology in graph filtration. The 0D and 1D persistent diagrams  $(b_i, d_i)$ tabulate the life-time of 0-cycles (connected components) and 1-cycles (loops) that are born at the filtration value $b_i$ and die at value $d_i$, respectively. The 0th Betti number $\beta_0(w_{(i)})$ counts the number of 0-cycles at filtration value $w_{(i)}$  and can be shown to be non-decreasing over filtration (Figure \ref{fig:BDschematic}) \citep{chung.2019.ISBI}: 
$$\beta_0(w_{(i)}) \leq \beta_0(w_{(i+1)}).$$
On the other hand the 1st Betti number $\beta_1(w_{(i)})$ counts the number of independent loops and can be shown to be non-increasing over filtration  \citep{chung.2019.ISBI}:
$$\beta_1(w_{(i)}) \geq \beta_1(w_{(i+1)}).$$
During the graph filtration, when new components is born, they never die. Thus, 0D persistent diagrams are completely characterized by birth values $b_i$ only. Loops are viewed as already born at $-\infty$. Thus, 1D persistent diagrams are completely characterized by death values $d_i$ only. We can show that the edge weight set $W$ can be partitioned into 0D birth values and 1D death values \citep{song.2021.MICCAI}:

\begin{figure}[t]
\begin{center}
\includegraphics[width=1\linewidth]{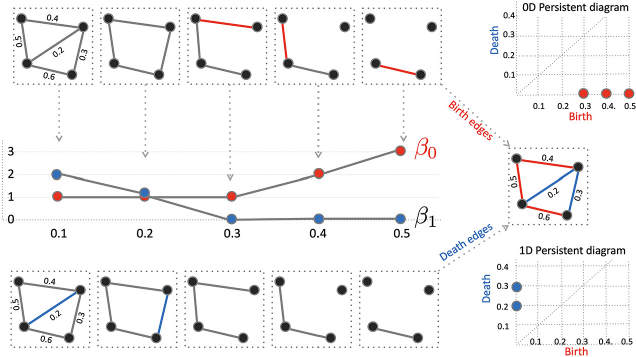}
\caption{The birth-death decomposition partitions the edge set into the birth and death edge sets. The birth set forms the maximum spanning tree (MST) and contains edges that create connected components (0D topology). The death set contains edges that do not belong to the maximum spanning tree (MST) and destroys loops (1D topology).}
\label{fig:BDschematic}
\end{center}
\end{figure}

\begin{theorem}[Birth-death decomposition]
The edge weight set $W =  \{ w_{(1)}, \cdots, w_{(q)} \}$ has the unique decomposition
\bqn W = W_b \cup W_d,  \quad  W_b \cap W_d = \emptyset, \label{eq:decompose} \eqn
where birth set $W_b = \{ b_{(1)}, b_{(2)}, \cdots, b_{(q_0)} \}$ is the collection of 0D sorted birth values and death set $W_d = \{ d_{(1)}, d_{(2)}, \cdots, d_{(q_1)} \}$ is the collection of 1D sorted death values with $q_0 = p-1$ and $q_1 = (p-1)(p-2)/2$. Further $W_b$ forms the 0D persistent diagram while $W_d$ forms the 1D persistent diagram. \label{thm:decompose}
\end{theorem}

In a complete graph with $p$ nodes, there are $q=p(p-1)/2$ unique edge weights. There are $q_0 = p-1$ number of edges that produce 0-cycles. This is equivalent to the number of edges in the maximum spanning tree (MST) of the graph. Thus, $$q_1 = q - q_0 = \frac{(p-1)(p-2)}{2}$$ number of edges destroy loops. The 0D persistent diagram  is given by $\{ (b_{(1)}, \infty),$ $\cdots,$ $(b_{(q_0)}, \infty) \}$. Ignoring $\infty$, $W_b$ is the 0D persistent diagram. The 1D persistent diagram  is given by $\{ (-\infty, d_{(1)}),$ $\cdots, (-\infty, d_{(q_1)}) \}$. Ignoring $-\infty$, $W_d$ is the 1D persistent digram. We can show that the birth set is the MST (Figure \ref{fig:BDschematic}) \citep{song.2023}. \\

\begin{figure}[t]
\begin{center}
\includegraphics[width=1\linewidth]{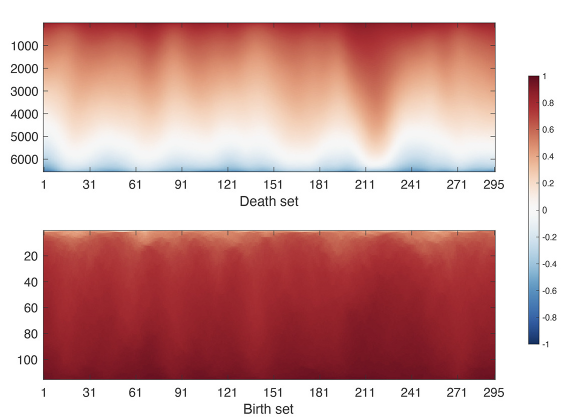}
\caption{The corresponding birth and death sets of dynamically changing correlation matrix shown in Figure \ref{fig:dynamicTDA}. The horizontal axis is the time point. Columns are the sorted birth and death edge values at the time point.}
\label{fig:BDtime}
\end{center}
\end{figure}

The identification of $W_b$ is based on the modification to Kruskal's or Prim's algorithm that identifies the MST \citep{lee.2012.tmi,song.2023}. Then $W_d$ is identified as $W \setminus W_d = W \cap W_d^c$. Figure  \ref{fig:BDtime} displays how the birth and death sets change over time for a single subject used in the study. Given edge weight matrix $W$ as input, the Matlab function {\tt WS\_decompose.m} outputs the birth set $W_b$ and the death set $W_d$.

\subsection{Topological distances}
\label{sec:dist}
\begin{figure}[t]
\begin{center}
\includegraphics[width=1\linewidth]{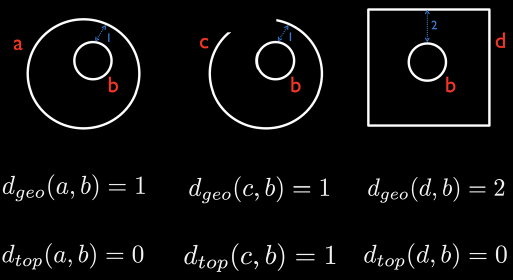}
\caption{Comparison between geometric distance $d_{geo}$ and topological distance $d_{top}$. We used the shortest Euclidean distance between objects as the geometric distance. The left (two circles) and middle (circle and arc) objects are topologically different while the left and right (square and circle) objects are topologically equivalent. The geometric distance cannot discriminate topologically different objects (left and middle) and produces false negatives. The geometric distance incorrectly discriminates topologically equivalent objects  (left and right) and produces false positive.}
\label{fig:geotop}
\end{center}
\end{figure}

Like the majority of clustering methods such as $k$-means and hierarchical clustering that use geometric distances \citep{johnson.1967,hartigan.1979,lee.2012.tmi}, we propose to develop a topological clustering method using topological distances (Figure \ref{fig:geotop}). For this purpose we use the Wasserstein distance. 

Given two probability distributions $X \sim f_1$ and $Y \sim f_2$, the $r$-{\em Wasserstein distance} $D_W$, the probabilistic version of the optimal transport,  is defined as
\bqn D_W (f_1, f_2) =  \Big( \inf \mathbb{E} |  X - Y |^r \Big)^{1/r}, \eqn
where the infimum is taken over every possible joint distribution of $X$ and $Y$.  The Wasserstein distance
is the optimal expected cost of transporting points generated from $f_1$ to those generated from $f_2$ \citep{canas.2012}. There are numerous distances and similarity measures defined between probability distributions such as the Kullback-Leibler (KL) divergence and the mutual information \citep{kullback.1951}. While the Wasserstein distance is a metric satisfying positive definiteness, symmetry,
and triangle inequality, the KL-divergence and the mutual information are not metrics. Although they are easy to compute, the biggest limitation of the KL-divergence and the mutual information is that the two probability distributions must be defined on the same sample space. If the two distributions do not have the same support, it may be difficult to even define the distance between them. If $f_1$ is discrete while $f_2$ is continuous, it is difficult to define them. On the other hand, the Wasserstein distance can be computed for any arbitrary distributions that may not have the common sample space, making it extremely versatile.

Consider persistent diagrams $P_1$ and $P_2$ given by  
$$P_1: x_1 = (b_1^1, d_1^1), \cdots, x_q = (b_q^1,d_q^1),  \quad P_2: y_1 = (b_1^2, d_1^2), \cdots, y_q = (b_q^2, d_q^2).$$ 
Their empirical distributions are given in terms of Dirac-Delta functions 
$$f_1 (x) = \frac{1}{q} \sum_{i=1}^q \delta (x-x_i), \quad f_2(y) = \frac{1}{q} \sum_{i=1}^q \delta (y-y_i).$$
Then we can show that the $r$-{\em Wasserstein distance} on persistent diagrams  is given by
\bqn D_W(P_1, P_2) = \inf_{\psi: P_1 \to P_2} \Big( \sum_{x \in P_1} \| x - \psi(x) \|^r \Big)^{1/r} \label{eq:D_Winf} \eqn
over every possible bijection $\psi$, which is a permutation, between $P_1$ and $P_2$ \citep{vallender.1974,canas.2012,berwald.2018}. 
Optimization (\ref{eq:D_Winf}) is the standard assignment problem, which is  usually solved by the Hungarian algorithm in $\mathcal{O} (q^3)$ \citep{edmonds.1972}. However, for graph filtration, the distance can be computed {\em exactly} in $\mathcal{O}(q \log q)$  by simply matching the order statistics on the birth or death values \citep{rabin.2011,song.2023,song.2021.MICCAI}:

\begin{theorem} The r-Wasserstein distance between the 0D persistent diagrams for graph filtration is given by 
\bqn D_{W0}(P_1, P_2) = \Big[ \sum_{i=1}^{q_0} (b_{(i)}^1 - b_{(i)}^2)^r \Big]^{1/r}, \eqn
where $b_{(i)}^j$ is the $i$-th smallest birth values in persistent diagram $P_j$.  The $r$-Wasserstein distance between the 1D persistent diagrams for graph filtration is given by 
\bqn D_{W1}(P_1, P_2) =  \Big[ \sum_{i=1}^{q_1} (d_{(i)}^1 - d_{(i)}^2)^r  \Big]^{1/r}, \eqn
where $d_{(i)}^j$ is the $i$-th smallest death values in persistent diagram $P_j$. 
\label{theorem:optimal}
\end{theorem}

The proof is provided in \cite{chung.2023.DE}. We can show that the 2-Wasserstein distance is equivalent to the Euclidean distance within a certain convex set.
Let ${\bf b}_i=(b_{(1)}^{i},  b_{(2)}^{i}  \cdots, d_{(q_0)}^{i})^{\top}$ be the vector of sorted birth values of persistent diagram $P_i$. Then ${\bf b}_i$ is a point in the $(q_0-1)$-simplex $\mathcal{T}_0$ given by
$$  \mathcal{T}_0 = \{ (x_1, x_2, \cdots, x_{q_0})  | x_1 < x_2 < \cdots < x_{q_0} \} \subset \mathbb{R}^{q_0},$$
where $x_1$ and $x_{q_0}$ are bounded below and above respectively. If brith and death values are from correlation matrices, 
$-1 \leq x_1$ and $x_{q_0} \leq 1$.
Hence, the 0D Wasserstein distance is equivalent to Euclidean distance in the $q_0$-dimensional convex set $\mathcal{T}_0$.
Similarly,  the vector of  sorted death values ${\bf d}_i=(d_{(1)}^{i},  d_{(2)}^{i}  \cdots, d_{(q_1)}^{i})^{\top}$ of persistent diagram $P_i$  is a point in the $(q_1-1)$-simplex  $\mathcal{T}_1$ given by
\bqn \mathcal{T}_1 = \{ (x_1, x_2, \cdots, x_{q_1})  | x_1 < x_2 < \cdots < x_{q_1} \} \subset \mathbb{R}^{q_1}, \eqn
where $x_1$ and $x_{q_1}$ are bounded below and above respectively. Hence, the 1D Wasserstein distance is equivalent to Euclidean distance in the $q_1$-dimensional convex set $\mathcal{T}_1$.

\subsection{Topological mean and variance}

Given a collection of graphs $\mathcal{X}_1= (V, w^1), \cdots, \mathcal{X}_n= (V, w^n)$ with edge weights $w^k = (w_{ij}^k)$, the usual approach for obtaining the average network  $\bar{\mathcal{X}}$ is simply averaging the edge weight matrices in an element-wise fashion
$$\bar{\mathcal{X}}  =   \Big( V, \frac{1}{n} \sum_{k=1}^n w_{ij}^k \Big).$$
However, such average is the average of the connectivity strength. Such an approach is usually sensitive to topological outliers \citep{chung.2019.ISBI}. We address the problem through the Wasserstein distance. A similar concept was proposed in the persistent homology literature through the Wasserstein barycenter \citep{agueh.2011,cuturi.2014}, which is motivated by the Fr\'{e}chet mean \citep{le.2000,turner.2014,zemel.2019,dubey.2019}. However, the method has not seen many applications in modeling graphs and networks.

To account for both 0D and 1D topological differences in networks, we use the sum of 0D and 1D Wasserstein distances between networks \(\mathcal{X}_1\) and \(\mathcal{X}_2\) as the topological distance
\begin{equation} 
\mathcal{D}(\mathcal{X}_1, \mathcal{X}_2) = D_{W0}^2(P_1, P_2)  + D_{W1}^2(P_1, P_2). 
\label{eq:Dx} 
\end{equation}
The equal weights of the form (\ref{eq:Dx}) were used for the following reasons. Through the birth-death decomposition, a weighted graph can be topologically characterized by 0D and 1D features, with no higher-dimensional features present. However, it is unclear which feature contributes the most. Equal weighting of 0D and 1D features ensures a balanced representation without bias towards either type of feature. 

Let \(\boldsymbol{\varnothing}\) denote a graph with zero edge weights. Then, due to the birth-death decomposition, we have
\[
\mathcal{D}(\mathcal{X}_1, \boldsymbol{\varnothing})  = \sum_{i<j} (w_{ij}^1)^2, \quad \mathcal{D}(\boldsymbol{\varnothing}, \mathcal{X}_2)  = \sum_{i<j} (w_{ij}^2)^2,
\]
where the squared sums of all the edge weights make the interpretation straightforward. If unequal weighting is used, these relationships do not hold.
Further, the Wasserstein distances \(D_{W0}\) and \(D_{W1}\) are equivalent to the Euclidean distances in a convex set. Therefore, squared distance is a more natural choice that satisfies the triangle inequality
\[
\mathcal{D}(\mathcal{X}_1, \mathcal{X}_3) \leq \mathcal{D}(\mathcal{X}_1, \mathcal{X}_2) + \mathcal{D}(\mathcal{X}_2, \mathcal{X}_3),
\]
thus qualifying as a metric.

The sum (\ref{eq:Dx}) does not uniquely define networks. Like the toy example in Figure \ref{fig:geotop}, we can have many topologically equivalent brain networks that give the identical distance. Thus, the average of two graphs is also not uniquely defined. The situation is analogous to Fr\'{e}chet mean, which frequently does not result in a unique mean\citep{le.2000,turner.2014,zemel.2019,dubey.2019}. We introduce the concept of the {\em topological mean} for networks, defined as the minimizer according to the Wasserstein distance, mirroring how the sample mean minimizes the Euclidean distance. The squared Wasserstein distance is translation invariant such that
$$\mathcal{D}(c + \mathcal{X}_1, c+ \mathcal{X}_2) = \mathcal{D}(\mathcal{X}_1, \mathcal{X}_2).$$
If we scale connectivity matrices by $c$, we have
$$\mathcal{D}(c \mathcal{X}_1, c\mathcal{X}_2) = c^2 \mathcal{D}(\mathcal{X}_1, \mathcal{X}_2).$$

\begin{definition} 
\label{definition:means}
The topological mean $\mathbb{E} \mathcal{X}$ of networks  $\mathcal{X}_1, \cdots, \mathcal{X}_n$  is 
 the graph given by
\bqn \mathbb{E} \mathcal{X}  =     \arg \min_{X} \sum_{k=1}^n   \mathcal{D}( X, \mathcal{X}_k).\label{eq:topology-mean}\eqn
\end{definition}
Unlike the sample mean, we can have many different networks with identical topology that give the minimum. Similarly, we can define the  {\em topological variance} $\mathbb{V} \mathcal{X}$ as follows.
\begin{definition} The topological variance $\mathbb{V} \mathcal{X}$ of networks  $\mathcal{X}_1, \cdots, \mathcal{X}_n$  is 
 the graph given by
\bqn \mathbb{V} \mathcal{X} = \frac{1}{n} \sum_{k=1}^n \mathcal{D}(\mathbb{E} \mathcal{X}, \mathcal{X}_k). \eqn
\end{definition} 
The topological variance can be interpreted as the variability of graphs from the topological mean $\mathbb{E} \mathcal{X}$. To compute the topological mean and variance, we only need to identify a network with identical topology as the topological mean or the topological variance.

\begin{theorem} 
Consider graphs $\mathcal{X}_i  = (V, w^i)$ with corresponding birth-death decompositions $W_i =W_{ib} \cup W_{id}$ with 
birth sets $W_{ib} = \{  b_{(1)}^i, \cdots, b_{(q_0)}^i\}$ and death sets $W_{id} = \{  d_{(1)}^i, \cdots, d_{(q_1)}^i\}$. Then, there exists topological mean $\mathbb {E} \mathcal{X} =(V,w)$ with birth-death decomposition $W_b \cup W_d$ with 
$W_b = \{  b_1, \cdots, b_{q_0}\}$ and  $W_d = \{  d_1, \cdots, d_{q_1} \}$ satisfying
\bqn b_j = \frac{1}{n} \sum_{i=1}^n b_{(j)}^i, \quad d_j = \frac{1}{n} \sum_{i=1}^n d_{(j)}^i. \eqn
\label{theorem:TM}
\end{theorem}

\subsection{Topological clustering}
\label{sec:clustering}
There are few studies that used the Wasserstein distance for clustering \citep{mi.2018,yang.2020}. 
The existing methods are mainly applied to geometric data without topological consideration or persistence. It is not obvious how to apply such geometric methods to  cluster graph or network data. We propose to use  the Wasserstein distance to cluster collection of graphs $\mathcal{X}_1, \cdots, \mathcal{X}_n$  into $k$ clusters 
$C_1, \cdots,C_k$ such that 
$$\cup_{i=1}^k C_i = \{\mathcal{X}_1, \cdots, \mathcal{X}_n\} , \quad C_i \cap C_j = \emptyset.$$

Let $C=(C_1, \cdots, C_k)$ be the collection of clusters. Let $\mu_j$ be the {\em topological cluster mean} within $C_j$ given by
 $$\mu_j = \arg \min_{X} \sum_{\mathcal{X}_k \in C_j}   \mathcal{D}( X, \mathcal{X}_k).$$
The cluster mean is computed through Theorem \ref{theorem:TM}. Just like Fr\'{e}chet mean, the cluster mean is not unique in a geometric sense but only unique in a topological sense \citep{turner.2014,le.2000,zemel.2019,dubey.2019}. Let $\mu = (\mu_1, \cdots, \mu_k)$ be the cluster mean vector. The within-cluster distance from the cluster centers is given by
\bqn l_W (C; \mu) =  \sum_{j=1}^k \sum_{X \in C_j} \mathcal{D}(X, \mu_j).\label{eq:l_W}\eqn 
If we let $|C_j|$ to be the number of networks within cluster $C_j$,  (\ref{eq:l_W}) can be written as
\bqn l_W (C; \mu) =  \sum_{j=1}^k |C_j| \mathbb{V}_{j} \mathcal{X},  \eqn 
with topological cluster variance 
$$\mathbb{V}_{j} \mathcal{X} = \frac{1}{|C_j|} \sum_{X \in C_j} \mathcal{D}(X, \mu_j)$$ within cluster $C_j$. 
The optimal cluster is found by minimizing within-cluster distance $l_W(C;\mu)$ in (\ref{eq:l_W}) over every possible partition of $C$.

If $\mu$ is  given and fixed, the identification of clusters $C$ can be done easily by assigning each network to the closest mean. Thus the topological clustering algorithm can be written as the two-step optimization similar to the expectation maximization (EM) algorithm often used in variational inferences and likelihood methods  \citep{bishop.2006}.  The first step computes the cluster mean. The second step minimizes the within-cluster distance. Just like $k$-means clustering, the two-step optimization is then iterated until convergence. Such process converges locally. 

\begin{theorem}
\label{th:converge}
The topological clustering converges locally.
\end{theorem}

The direct algebraic proof is fairly involving and given in \citet{chung.2023.NI}. Here we provide a more intuitive explanation. Note $D_{W0}$ and $D_{W1}$ are Euclidean distances in  convex set  $\mathcal{T}_0 \otimes \mathcal{T}_1$. Subsequently, 
$$\mathcal{D}(P_1,P_2 ) = D_{W0}^2(P_1, P_2)+  D_{W1}^2(P_1,P_2).$$
is the Euclidean distance in the Cartesian product $\mathcal{T}_0 \otimes \mathcal{T}_1$.
Thus, our topological clustering is equivalent to $k$-means clustering restricted to the convex set $\mathcal{T}_0 \otimes \mathcal{T}_1$. The convergence of topological clustering is then the direct consequence of the convergence of $k$-means clustering, which always converges in such a convex space. Numerically we minimize (\ref{eq:l_W}) by replacing the 
Wasserstein distance with the 2-norm between sorted vectors of birth and death values in $k$-means clustering.

Like $k$-means clustering algorithm that only converges to local minimum, there is no guarantee the topological clustering converges to the global minimum \citep{huang.2020.NM}. This is remedied by repeating the algorithm multiple times with different random seeds and taking the smallest possible minimum. The method is implemented as the Matlab function {\tt WS\_cluster.m} which inputs the collection of networks and  outputs the cluster labels and clustering accuracy.  Different choice of initial cluster centers may lead to different results. Thus, the algorithm may become stuck in a local minimum and may not converge to the global minimum. Thus, in actual numerical implementation, we used different initializations of centers. Then, we picked the best  clustering result with the smallest within cluster distance $l_W$.

\subsection{Validation}

\begin{figure}[t]
\centering
\includegraphics[width=0.8\linewidth]{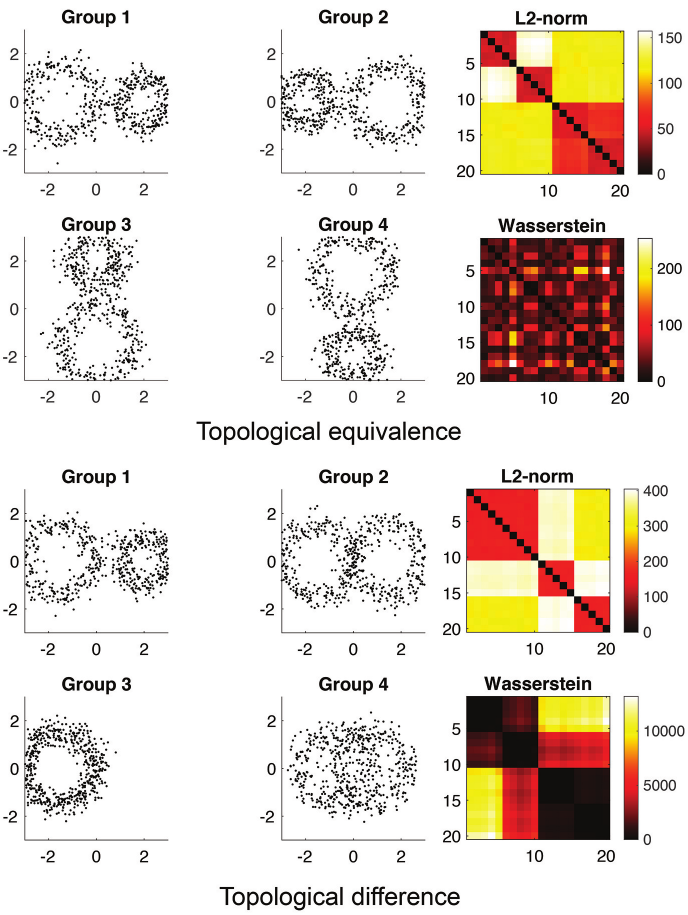}
\caption{Top: simulation study on  topological equivalence. The correct clustering method should {\em not} be able to cluster them because they are all topologically equivalent. The pairwise Euclidean distance ($L_2$-norm) is used in $k$-means and hierarchical clustering. The Wasserstein distance is used in topological clustering. Bottom: simulation study on  topological difference. The correct clustering method should be able to cluster them because they are all topologically different.}
\label{fig:simulationnodiff}
\end{figure}

We validated the topological clustering in a simulation with the ground truth against $k$-means  and hierarchical clustering \citep{lee.2011.MICCAI}. We generated 4 circular patterns of identical topology (Figure \ref{fig:simulationnodiff}-top) and different topology (Figure \ref{fig:simulationnodiff}-bottom). Along the circles, we uniformly sampled 60 nodes and added  Gaussian noise $N(0, 0.3^2)$ on the coordinates. We generated 5 random networks per group. The Euclidean distance ($L_2$-norm)  between randomly generated points are used to build connectivity matrices for $k$-means  and hierarchical clustering. Figures \ref{fig:simulationnodiff} shows the superposition of nodes from 20 networks.   For $k$-means and Wasserstein graph clustering, the average result of 100 random seeds is reported.

  \subsubsection{Testing for false positives}

In the experiment depicted in Figure \ref{fig:simulationnodiff}, we evaluated the occurrence of false positives in scenarios devoid of topological differences. All groups, derived from Group 1 through rotations, are topologically identical. Hence, any detected differences are false positives. 
While $k$-means clustering exhibited an accuracy of $0.90 \pm 0.15$, and hierarchical clustering achieved perfect accuracy (1.00),
these methods reported significant false positives, erroneously categorizing the groups as distinct clusters. The absence of inherent topological differences between the groups implies that higher clustering accuracy is indicative of false positive results.
Conversely, topological clustering, with a lower accuracy of $0.53 \pm 0.08$, demonstrated a reduced tendency for reporting false positives  in the absence of topological differences.

\subsubsection{Testing for false negatives}

Figure \ref{fig:simulationnodiff}  presents our test for false negatives, featuring groups with varying numbers of cycles and distinct topologies. In this scenario, topological differences should be detectable. Here, $k$-means clustering recorded an accuracy of $0.83 \pm 0.16$, and hierarchical clustering again reported perfect accuracy. Notably, topological clustering attained a high accuracy of $0.98 \pm 0.09$. Separating topological from geometric signals is challenging; the presence of topological differences often coincides with geometric variations, which can influence the performance of all tested methods.

In summary, while traditional clustering methods based on geometric distances are prone to a significant number of false positives, making them less suitable for topological learning tasks, the proposed Wasserstein distance-based approach demonstrates superior performance. This method excels in minimizing both false positives and false negatives, as evidenced by our tests. Its effectiveness is particularly noteworthy in topological learning tasks, where discerning topological rather than geometric distinctions is crucial.

\subsection{Weighted Fourier series representation}

The predominant method for computing time-varying correlation in time series data, particularly in neuroimaging studies, involves Sliding Windows (SW). This technique entails computing correlations between brain regions across various time windows \citep{allen.2014, hutchison.2013, shakil.2016, mokhtari.2019, huang.2020.NM}. However, the use of discrete windows in SW can lead to artificially high-frequency fluctuations in dynamic correlations \citep{oppenheim.2001}. While tapering methods can occasionally mitigate these effects \citep{allen.2014}, the correlation computations within these windows remain susceptible to the influence of outliers \citep{devlin.1975}.

To circumvent these limitations, we employed the Weighted Fourier Series (WFS) representation \citep{chung.2007.TMI,chung.2008.TMI}. This approach extends the traditional cosine Fourier transform by incorporating an additional exponential weight. This modification effectively smooths out high-frequency noise and diminishes the Gibbs phenomenon \citep{chung.2007.TMI, huang.2019.DSW}. Crucially, WFS eliminates the need for sliding windows (SW) when computing time-correlated data. Given the necessity for robust signal denoising methods to ensure the efficacy of the persistent homology method across various subjects and time points, such an approach is needed. Consider an arbitrary noise signal \( f(t), t \in [0,1] \), which will undergo denoising through the diffusion process.

\begin{theorem} The unique solution to 1D heat diffusion:
\bqn \frac{\partial}{\partial {s}} h(t,s) = \frac{\partial^2}{\partial t^2} h(t,s) \label{eq:diffusion1} \eqn
on unit interval $[0,1]$ with initial condition $h(t,{s}=0)=f(t)$ is given by  WFS:
\bqn h(t,{s})=\sum_{l=0}^\infty e^{- l^2\pi^2 {s}}c_{fl}\psi_l(t), 
\label{eq:hts01} 
\eqn
where $\psi_0(t) =1$, $\psi_l(t) = \sqrt{2} \cos ( l \pi t)$ are the cosine basis 
and  $ c_{fl}=\int_{0}^1f(t)\psi_l(t)dt$ are the expansion coefficients.
\end{theorem}

The algebraic derivation is given in \citep{chung.2007.TMI}. Note the cosine basis is orthonormal
$$\langle \psi_l, \psi_m \rangle = \int_0^1 \psi_l(t) \psi_m(t) \; dt = \delta_{lm},$$
where $\delta_{lm}$ is Kroneker-detal taking value 1 if $l=m$ and 0 otherwise.  We can rewrite (\ref{eq:hts01}) as a more convenient convolution form
\bq h(t,{s})=\int_{0}^1K_s(t,t')f(t')dt', \label{eq:hts}\eq
where heat kernel $K_s(t,t')$ is given by
\bq
K_s(t,t')=\sum_{l=0}^\infty e^{- l^2\pi^2 {s}}\psi_l(t)\psi_l(t').
\label{eq:heat}
\eq
The diffusion time $s$ is usually referred to as the kernel bandwidth and  controls the amount of smoothing. 
Heat kernel satisfies $\int_0^1 K_s(t,t') \; dt = 1$ for any $t'$ and $s$.

To reduce unwanted boundary effects near the data boundary $t=0$ and $t=1$ \citep{huang.2019.DSW,huang.2020.NM}, we project the data onto  the circle $\mathcal{C}$ with circumference 2 by the mirror reflection:
$$g(t) = f(t) \mbox{ if } t \in [0,1], \quad g(t)=f(2-t) \mbox{ if } t \in [1,2].$$ 
Then perform WFS on the circle.

\begin{figure}[t]
\centering
\includegraphics[width=1\linewidth]{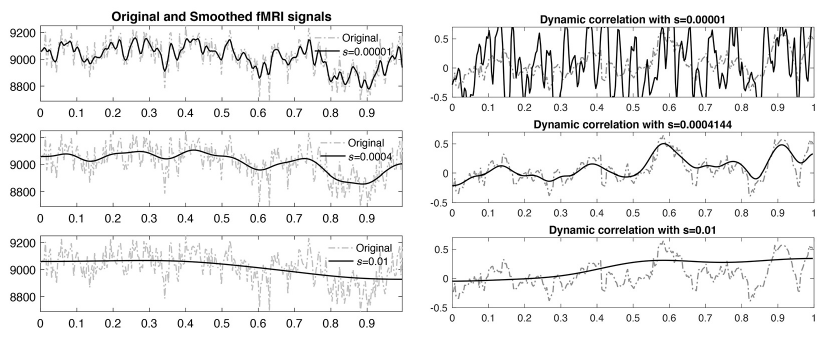}
\caption{Left: The original and smoothed fMRI time series using WFS with degree $L = 295$ and different heat kernel bandwidth $s$. The bandwidth $4.141 \times 10^{-4}$ used in this study approximately matches 20 TRs often used in the sliding window methods. Right: The dotted gray lines are correlations computed over sliding windows. The solid black lines are correlations computed  using WFS.}
\label{fig:HeatCorr}
\end{figure}

\begin{theorem} 
\label{th:diffcircle}
The unique solution to 1D heat diffusion:
\bqn \frac{\partial}{\partial {s}} h(t,s) = \frac{\partial^2}{\partial t^2} h(t,s) \label{eq:diffusion2} \eqn
on the circle $\mathcal{C}$  with the initial periodic condition
$h(t,s=0)  = f(t) \mbox{ if } t \in [0,1], h(t,s=0) =f(2-t) \mbox{ if } t \in [1,2]$ is given by  WFS:
\bqn h(t,{s})= \sum_{l=0}^\infty e^{- l^2\pi^2 {s}}c_{fl}\psi_l(t),
\eqn
where $\psi_0(t) =1$, $\psi_l(t) = \sqrt{2} \cos ( l \pi t)$ are the cosine basis 
and  $ c_{fl}=\int_{0}^1f(t)\psi_l(t)dt$ are the expansion coefficients.
\end{theorem}

The cosine series coefficients $c_{fl}$ are  estimated using the least squares method by setting up a matrix equation \citep{chung.2007.TMI}. We set the expansion degree to equate the number of time points, which is 295. The window size of 20 TRs was used in most sliding window methods \citep{allen.2014,lindquist.2014,huang.2020.NM}. We matched the full width at half maximum (FWHM) of heat kernel to the window size numerically. We used the fact that diffusion time $s$ in heat kernel approximately matches to the kernel bandwidth  of Gaussian kernel $e^{-t^2/2\sigma^2}$ as $\sigma = s^2/2$ (page 144 in \citep{chung.2012.CNA}). 20 TRs  is approximately equivalent to heat kernel bandwidth of about $4.144 \cdot 10^{-4}$ in terms of FWHM. Figure \ref{fig:HeatCorr} displays the WFS representation of rsfMRI with different kernel bandwidths.

 \subsection{Dynamic correlation on weighted Fourier series}
 The weighted Fourier series representation provides a way to compute correlations dynamically without using sliding windows. Consider time series $x(t)$ and $y(t)$ with heat kernel $K_s(t,t')$. The mean and variance of signals with respect to the heat kernel are given by
\bq \mathbb{E} x(t) &=& \int_0^1  K_s(t,t') x(t') \; dt',\\
\mathbb{V} x(t) &=& \int_0^1  K_s(t,t') x^2(t') \; dt' - \big[ \mathbb{E} x(t) \big]^2.
\eq
Subsequently, the correlation $w(t)$ of $x(t)$ and $y(t)$ is given by 
\bqn w(t) = \frac{\int_0^1 K_s(t,t') x(t')y(t') \; dt - \mathbb{E} x(t) \mathbb{E} y(t) }
{ \sqrt{ \mathbb{V} x(t)} \sqrt{ \mathbb{V} y(t)}}.
\eqn
When the kernel is shaped as a sliding window, the correlation $w(t)$ exactly matches the correlation computed over the sliding window. The kernelized correlation generalizes the concept of integral correlations with the additional weighting term \citep{huang.2019.ISBI}. As $s \to \infty$, $w(t)$ converges to the Pearson correlation computed over the whole time points. Thus, the kernel bandwidth behaves like the length of sliding window.

\begin{theorem}The correlation $w(t)$ of time series $x(t)$ and $y(t)$ with respect to heat kernel $K_s(t,t')$ is given by
\bqn \label{eq:windowless}
w (t) = \frac{\sum_{l=0}^\infty e^{-l^2 \pi^2 s}  c_{xyl}\psi_l(t) -\mu_x(t)\mu_y(t)}
{\sigma_{x}(t) \sigma_{y}(t) },
\eqn
with
\bq \mu_x(t) = \sum_{l=0}^\infty e^{-l^2 \pi^2 s} c_{xl}\psi_l(t), \quad \sigma_{x}^2(t) =\sum_{l=0}^\infty e^{-l^2 \pi^2 s}c_{xxl}\psi_l(t)-\mu^2_x(t).\eq
$$c_{xl}=\int_{0}^1x(t)\psi_l(t)dt, \quad c_{yl}=\int_{0}^1y(t)\psi_l(t)dt$$ are the cosine series coefficients. 
Similarly we expand $x(t)y(t)$, $x^2(t)$ and $y^2(t)$ using the cosine basis and obtain coefficients $c_{xyl}$, $c_{xxl}$ and $c_{yyl}$. 
\end{theorem}
The derivation follows by simply replacing all the terms with the WFS representation. Correlation (\ref{eq:windowless}) is the formula we used to compute the dynamic correlation in this study. Figure \ref{fig:HeatCorr} displays the WFS-based dynamic correlation for different bandwidths. A similar weighted correlation was proposed in  \citet{pozzi.2012}, where the time varying exponential weights proportional to $e^{t/\theta}$ with exponential decay factor $\theta$ were used. However, our exponential weight term is related to the spectral decomposition of heat kernel in the spectral domain and invariant over time. The WFS based correlation is not related to \citet{pozzi.2012}.

\begin{figure}[t]
\centering
\includegraphics[width=1\linewidth]{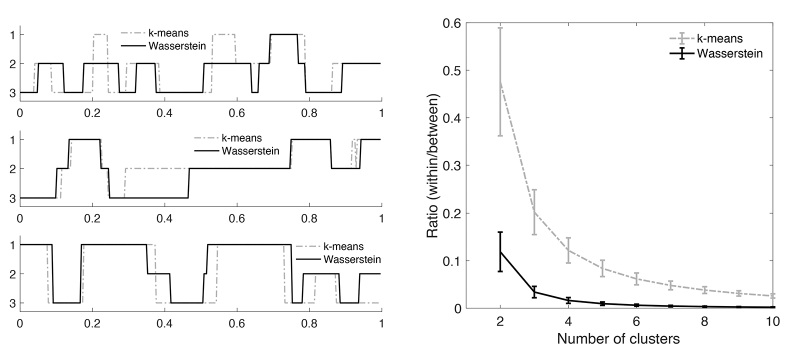}
\caption{Left:  The time series of estimated state spaces using the topological clustering and $k$-means clustering for 3 subjects. The time is normalized into unit interval $[0,1]$. Right: The ratio of within-cluster to between-cluster distances. Smaller the ratio, better the clustering fit is.}
\label{fig:ratio}
\end{figure}

\section{Results}

\subsection{Data}
The proposed method is applied in the accurate estimation of state spaces in dynamically changing functional brain networks. The 479 subjects resting-state functional magnetic resonance images (rs-fMRI) used in this paper were collected on a 3T
MRI scanner (Discovery MR750, General Electric Medical Systems, Milwaukee, WI, USA) with  a 32-channel RF head coil array. The 479 healthy subjects consist of 231 males and 248 females ranging in age from 13 to 25 years. The sample  contains 132 monozygotic (MZ) twin pairs and 93 same-sex dizygotic (DZ) twin pairs. 
 
 The image preprocessing includes motion corrections and image alignment to the MNI template and  follows   \citep{burghy.2016,jenkinson.2002}. The resulting rs-fMRI consist of  $91\times109\times91$  isotropic voxels at 295 time points. We further parcellated the brain volume into 116 non-overlapping brain regions using the Automated Anatomical Labelling (AAL) atlas \citep{tzourio.2002}. The fMRI data were averaged across voxels within each brain region, resulting in 116 average fMRI signals with 295 time points for each subject. The rs-fMRI signals were then scaled to fit to unit interval [0, 1] and treated as functional data in $[0,1]$.

\subsection{Topological state space estimation}\label{subsec:states}
For $p$ brain regions, we estimated $p\times p$ dynamically changing correlation matrices $C_i(t)$  for the $i$-th subject  using WFS. Let $\mathbf{C}_{ij}$ denote the vectorization of  the upper triangle of $p \times p$ matrix $C_i(t_j)$ at time point $t_j$ into $p^2 \times 1$ vector. For each fixed $i$, the collection of $\mathbf{C}_{ij}$ over $T$ = 295 time points  is then feed into topological clustering in identifying the recurring brain connectivity at the subject level. We are clustering individual brain networks without putting any constraint on group or twin. We compared  the proposed Wasserstein clustering against the $k$-means clustering, which has been often used as the baseline method in the state space modeling \citep{allen.2014,huang.2019.DSW,huang.2020.NM}. After clustering, each correlation matrix $C_i(t_j)$ is assigned integers between 1 and $k$. These discrete states serve as the basis for investigating  the dynamic pattern of brain connectivity \citep{ting.2018}. For the convergence of both topological clustering and $k$-means clustering, the clusterings were repeated 10 times with different initial centroids and  the best result (smallest within-cluster distance) is reported. Figure \ref{fig:ratio}-left displays the result of the topological clustering against the $k$-means for three subjects.  295 time points are rescaled to fit into unit interval $[0,1]$.

\begin{figure}[t]
\includegraphics[width=\linewidth]{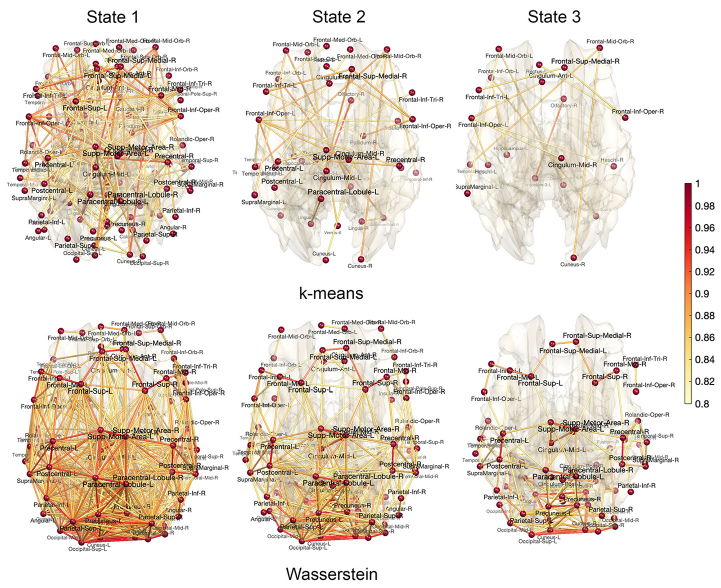}
\caption{The estimated state spaces of dynamically changing brain networks. The correlations are averaged over every time point and subject within each state for $k$-means clustering (top) and Wasserstein distance based topological clustering (bottom). In $k$-means clustering, the connectivity pattern of each state is somewhat random. 
In topological clustering,  the connectivity pattern is highly symmetric even though we did not put any symmetry constraint in the clustering method.}
\label{fig:clustering-results}
\end{figure}

The optimal number of cluster $k$ was determined by the {\em elbow method} \citep{allen.2014,rashid.2014,ting.2018,huang.2020.NM}. For each value of $k$, we computed the ratio of the within-cluster  to between-cluster distances. The ratio shows the goodness-of-fit of the cluster model.  
The elbow method  gives the largest slope change in the ratio when $k=3$  in the both methods (Figure \ref{fig:ratio}-right). At $k=3$, the ratio  is $0.034 \pm 0.012$ for 479 subjects for Wasserstein while it is $0.202 \pm   0.047$ for the $k$-means. The {\em six} times smaller ratio for the topological clustering  demonstrates the superior model fit over $k$-means. Figure \ref{fig:clustering-results} shows the results of clustering for both methods. 
 The $k$-means clustering result is based on averaging correlations of every time point and subject within each state.  The resulting states in the $k$-means clustering are somewhat random without any biologically interpretable pattern. The topological clustering computes the {\em topological mean} of every time point and subject within each state.

\subsection{Twin correlations over transpositions}

\begin{figure}[t]
\centering
\includegraphics[width=1\linewidth]{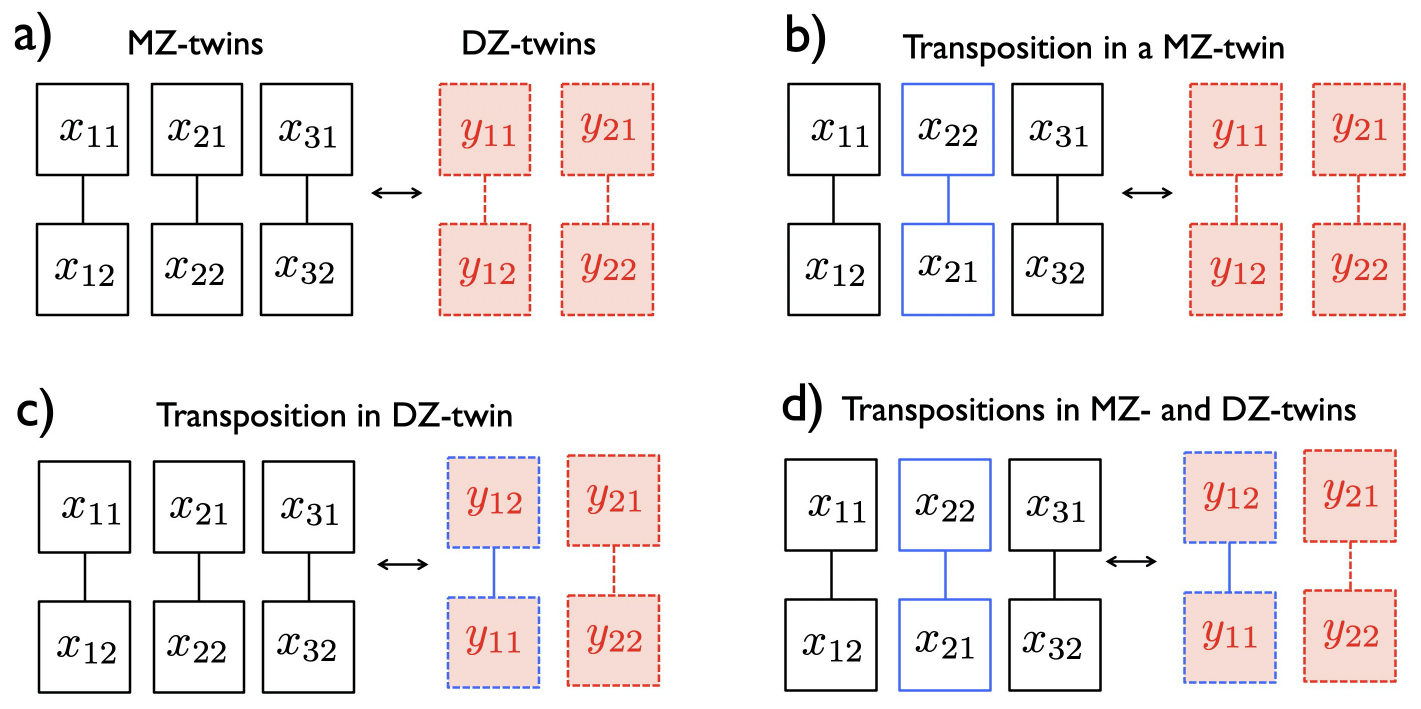}
\caption{The schematic of transpositions on 3 MZ- and 2 DZ-twins. a) One possible pairing. b) Transposition  within a MZ-twin. c) Transposition within a DZ-twin. d) Transpositions in both MZ- and DZ-twins simultaneously. Any transposition will affect the heritability estimate so it is necessary to account for as many transpositions as possible.}
\label{fig:schematic}
\end{figure}

Using additional twin information in the data, we further investigated if the state change pattern itself is genetically heritable. As far as we are aware, there is no study on the heritability of the state change pattern itself. This requires computing twin correlations. 
We assume there are $m$ MZ- and $n$ DZ-twins. For some feature, let $x_{i} = (x_{i1}, x_{i2})^{\top}$ be the  $i$-th twin pair in MZ-twin and $y_{i} = (y_{i1}, y_{i2})^{\top}$ be the  $i$-th twin pair in DZ-twin. They are represented as
\bq {\bf x} =
\left(
\begin{array}{ccc}
 x_{11}, & \cdots   &, x_{m1}   \\
 x_{12}, &  \cdots &,   x_{m2} \\
\end{array}
\right), \quad
{\bf y} =
\left(
\begin{array}{ccc}
 y_{11}, & \cdots   &, y_{n1}   \\
 y_{12}, &  \cdots &,   y_{n2} \\
\end{array}
\right).
\eq
Let ${\bf x}_j$ be the $j$-th row of ${\bf x}$, i.e., ${\bf x}_j = (x_{1j}, x_{2j}, \cdots, x_{mj})$. Similarly let  ${\bf y}_j = (y_{1j}, y_{2j}, \cdots, y_{nj})$. Then MZ- and DZ-correlations  are computed as the sample correlation
\bq \gamma^{MZ} ({\bf x}_1,{\bf x}_2) &=& corr({\bf x}_1,{\bf x}_2)\\
\gamma^{DZ} ({\bf y}_1,{\bf y}_2) &=& corr({\bf y}_1,{\bf y}_2).
\eq
However, there is no preference for the order of twins within a pair, and we can {\em transpose} the $i$-th twin pair in MZ-twin such that
\bq \tau_{i}({\bf x}_1) &=& (x_{11} \cdots x_{i-1,1}, x_{i2}, x_{i+1,1} \cdots x_{m1}),\\
 \tau_i({\bf x}_2) &=& (x_{12} \cdots x_{i-1,2}, x_{i1}, x_{i+1,2} \cdots x_{m2})
 \eq
and obtain another twin correlation $\gamma^{MZ}(\tau_i({\bf x}_1), \tau_i({\bf x}_2))$ \citep{chen.2018,chung.2019.CNI}.  Ignoring symmetry,  there are $2^{m}$ possible combinations in ordering the twins, which form a permutation group. The size of the permutation group grows exponentially large as the sample size increases.  Computing correlations over all permutations is not even computationally feasible for large $m$ beyond 100. Figure \ref{fig:schematic} illustrates many possible transpositions within twins. Thus, we propose a new fast online computational strategy for computing twin correlations.

Over transposition $\tau_i$, the correlation changes 
\bqn \gamma^{MZ}({\bf x}_1, {\bf x}_2) \to \gamma^{MZ}( \tau_i ({\bf x}_1), \tau_i ({\bf x}_2))
\eqn incrementally. We will determine the exact increment over the transposition. The sample correlation between ${\bf x}_k$ and ${\bf x}_l$ involves the following functions. 
\bq \nu ({\bf x}_k) &=&  \sum_{l=1}^m x_{lk}\\
  \omega({\bf x}_k, {\bf x}_l) &=&  \sum_{r=1}^m \big(x_{rk}-\nu({\bf x}_k)/m \big) \big(x_{rl}-\nu({\bf x}_l)/m \big). \eq
The functions $\mu$ and $\omega$ are updated over transposition $\tau_i$ as
\bq \nu (\tau_i({\bf x}_k)) &=& \nu({\bf x}_k) - x_{ik} + x_{il} \\
\omega(\tau_i({\bf x}_k), \tau_i ({\bf x}_l))&= &\omega({\bf x}_k, {\bf x}_l)+ (x_{ik} - x_{il})^2/m
 - (x_{ik} - x_{il}   )   \big(  \nu({\bf x}_k)  -\nu ({\bf x}_l) \big) /m.
\eq
Then the MZ-twin correlation after transposition is  updated as
\bqn
\gamma^{MZ}(\tau_i({\bf x}_1), \tau_i({\bf x}_2))= \frac{\omega(\tau_i({\bf x}_1), \tau_i({\bf x}_2))}{\sqrt{\omega(\tau_i({\bf x}_1), \tau_i({\bf x}_1)) \omega(\tau_i({\bf x}_2), \tau_i({\bf x}_2))}}. \label{eq:online-corr}
\eqn
The time complexity for correlation computation is 33 operations per transposition, which is substantially lower than the computational complexity of directly computing correlations per permutation. In the numerical implementation, we  sequentially apply random transpositions $\tau_{i_1}, \tau_{i_2}, \cdots, \tau_{i_J}$. This results in $J$ different twin correlations, which are averaged. Let $$\pi_1 = \tau_{i_1}, \pi_2 = \tau_{i_2} \circ \tau_{i_1}, \cdots, \pi_J = \tau_{i_J} \circ \cdots \circ  \tau_{i_2} \circ \tau_{i_1}.$$ 
The average correlation $\overline \gamma_J^{MZ}$ of all  $J$ transpositions is given by
$$ \overline \gamma_J^{MZ}= \frac{1}{J} \sum_{j=1}^J \gamma^{MZ} (\pi_{i_j}({\bf x}_1), \pi_{i_j}({\bf x}_2)).$$
In each sequential update, the average correlation can be updated iteratively as 
$$\overline \gamma_J^{MZ} =  \frac{J-1}{J} \overline \gamma_{J-1}^{MZ}  +  \frac{1}{J}\gamma^{MZ} (\pi_{i_J}({\bf x}_1), \pi_{i_J}({\bf x}_2)).$$
If we use enough transpositions, the average correlation $\overline \gamma_J^{MZ}$ converges to the true underlying twin correlation $\gamma^{MZ}$ for sufficiently large $J$. DZ-twin correlation $\gamma^{DZ}$ is estimated similarly.

In  the widely used ACE genetic model, the heritability index (HI) $h$, which  determines the amount of variation due to genetic difference in a population, is  estimated using Falconer's formula \citep{falconer.1995,chung.2019.NN,arbet.2020}. MZ-twins share 100\% of genes while same-sex DZ-twins share 50\% of genes on average. Thus, the additive genetic factor $A$ and the common environmental factor $C$  are related as 
\bq \gamma^{MZ} &=&A + C,\\ 
\gamma^{DZ} &=& A/2 + C.\label{eq:HI} 
\eq
HI $h$, which measures the contribution of $A$, is given by 
$$h({\bf x}, {\bf y}) = 2 ( \gamma^{MZ} - \gamma^{DZ} ).$$
In numerical implementation, 100 million transpositions can be easily done in 100 seconds in a desktop.  Similarly, we update the DZ-correlation over the transposition. 

\subsection{Heritability of the state space}

\begin{figure}[t]
\centering
\includegraphics[width=1\linewidth]{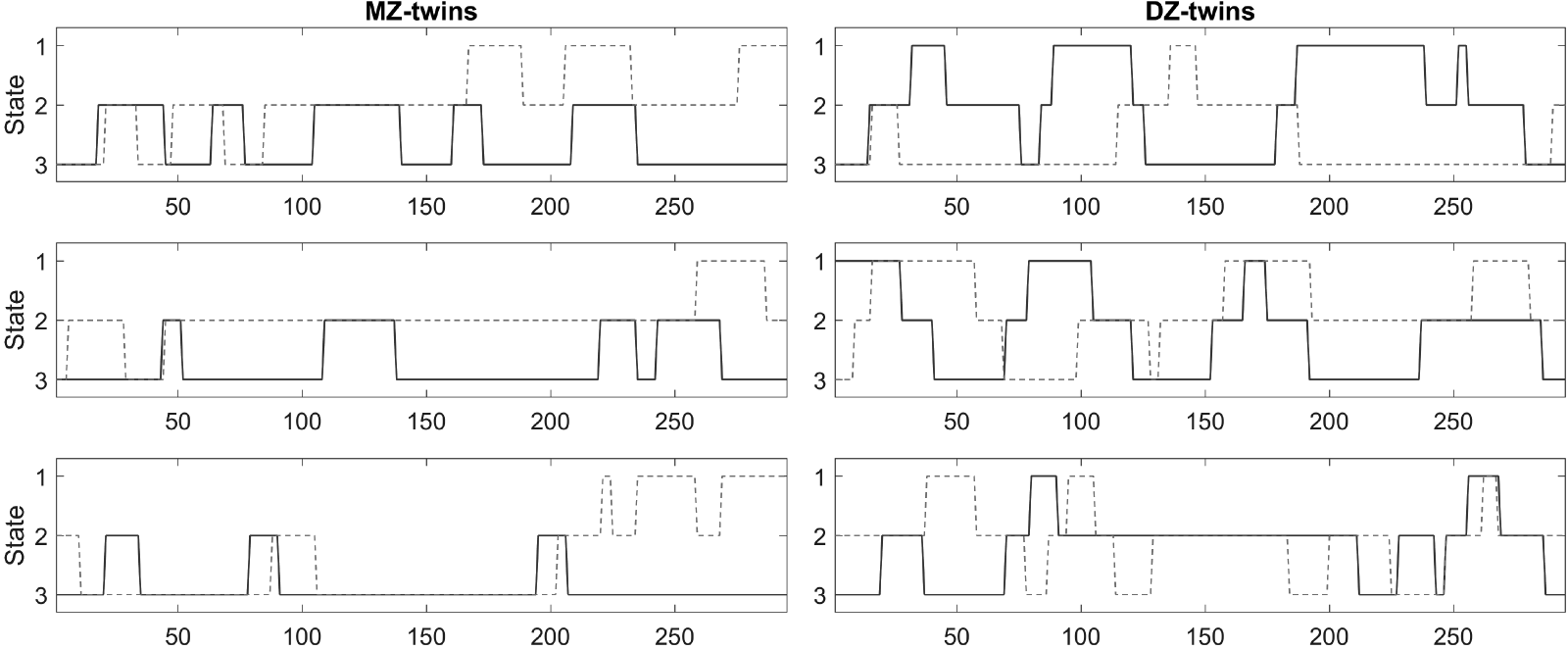}
\caption{State visits for 3 MZ-twins (left) and 3 DZ-twins (right) obtained from  the baseline $k$-means clustering. We are interested in determining the heritability of such state changes. Unfortunately, even within a twin, the time series of state change  do not synchronize, making the task extremely challenging.}
\label{fig:statevisits}
\end{figure}

The heritability estimation of state space is not a trivial task since the estimated state does not synchronize across twins making the task fairly difficult. Figure \ref{fig:statevisits} displays the state visits in randomly selected 3 MZ- and 3 DZ-twins. However, the time series of state changes do not synchronize within twins. This is likely a reason for the lack of reported  heritability of the state space in the literature. 

\begin{figure}[t]
\includegraphics[width=\linewidth]{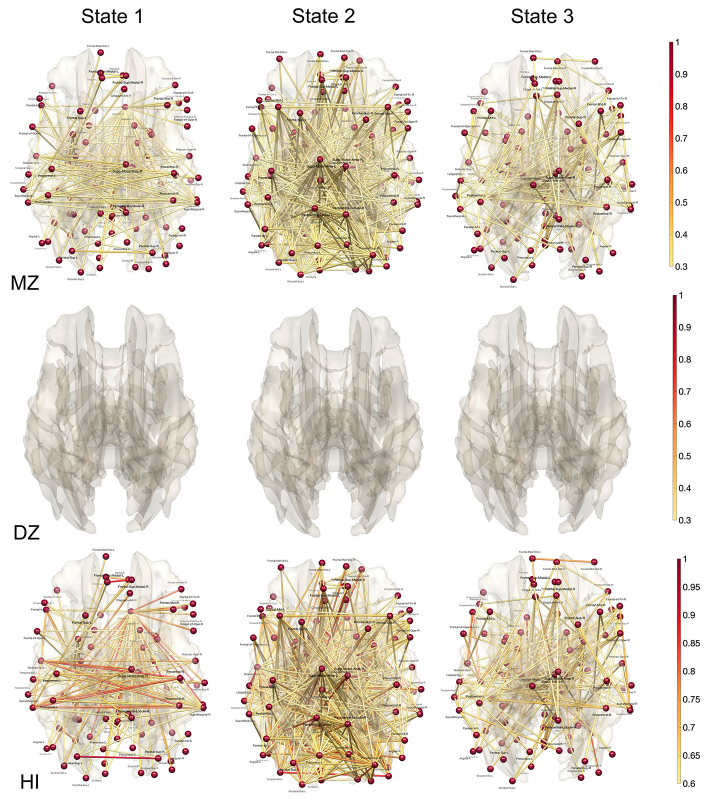}
\caption{\footnotesize MZ-correlation (top) and DZ-correlation (middle) in each state obtained through topological clustering in Figure \ref{fig:clustering-results}. There is no DZ-correlation above 0.3 and not displayed. The heritability index (HI) is determined by the twice the difference in twin correlations. HI of each state shows the extensive genetic contribution of dynamically changing states. The first state is characterized by prominent bilateral connections between the left and right hemispheres, whereas the second state primarily features front-back connections.}
\label{fig:HI}
\end{figure}

For each subject,  we computed the average correlation of each state, where the average is taken over all time points within each state. The correlation matrices are then used as the input to the transposition based twin correlations \citep{chung.2019.NN}. Subsequently, we computed the MZ- and DZ-twin correlations within each state (Figure \ref{fig:HI}). 
The MZ-twin correlations (Figure \ref{fig:HI}-top) are densely observed in many connections while there is no DZ-twin correlations (Figure \ref{fig:HI}-middle) observed above 0.3. We then computed the heritability index (HI) of each state (Figure \ref{fig:HI}-bottom). The heritability of the first state is characterized by strong lateralization of the hemisphere connections. The heritability of the second state is characterized by front and back connections. We believe the topological approach provides far more accurate and stable heritability index maps for dynamically changing state, which are biologically interpretable. 

We reported 10 connections that give the highest HI values in all three states in Tables \ref{tab:my_table1}, \ref{tab:my_table2} and \ref{tab:my_table3}.  Although numerous studies report high heritability for anatomical features such as gray matter density, there are few rs-fMRI studies reporting heritability of rs-fMRI \citep{glahn.2010,korgaonkar.2014}. Most of these studies report low HI compared to our high HI. \citep{glahn.2010} reported HI of 0.104 in the left cerebellum, 0.304 in the right cerebellum, 0.331 in the left temporal parietal region, 0.420 in the right  temporal parietal region. \citep{korgaonkar.2014} reported HI of 0.41 in the connection between the posterior cingulate cortex and right inferior parietal cortex in the default mode network  involving 79 MZ- and 46 same-sex DZ-twins. Other connections are all reporting very low HI below 0.24. We believe our topological method is clustering topologically similar functional network patterns and significantly boost genetic signals. 

\begin{table}[h]
\centering
\begin{tabular}{|l|l|c|}
\hline
Regions & Regions & HI \\
\hline
Parietal-Sup-L & Parietal-Sup-R & 0.96 \\
Frontal-Sup-Medial-L & Frontal-Sup-Medial-R & 0.90 \\
Olfactory-R & Temporal-Mid-R & 0.89 \\
Precentral-R & Rolandic-Oper-L & 0.88 \\
Olfactory-R & Temporal-Inf-R & 0.88 \\
Olfactory-R & Fusiform-R & 0.87 \\
Olfactory-R & Cerebelum-4-5-L & 0.87 \\
Precentral-R & SupraMarginal-L & 0.85 \\
Rolandic-Oper-L & Postcentral-R & 0.84 \\
Olfactory-R & Lingual-L & 0.84 \\
\hline
\end{tabular}
\caption{10 connections with the highest heritability index for state 1. Connections are sorted with respect to HI values.}
\label{tab:my_table1}
\begin{tabular}{|l|l|c|}
\hline
Regions & Regions &HI \\
\hline
Hippocampus-L & Cerebelum-4-5-L & 1.00\\
Olfactory-L & Fusiform-R & 0.92 \\
Precuneus-R & Cerebelum-Crus2-L & 0.90 \\
Occipital-Sup-L & Fusiform-L & 0.89 \\
Supp-Motor-Area-L & Cerebelum-Crus2-L & 0.88 \\
Occipital-Mid-L & Occipital-Mid-R & 0.87 \\
Thalamus-L & Cerebelum-9-L & 0.86 \\
Rolandic-Oper-L & Temporal-Sup-L & 0.85 \\
Paracentral-Lobule-L & Cerebelum-Crus2-L & 0.85 \\
Caudate-R & Cerebelum-Crus2-L & 0.85 \\
\hline
\end{tabular}
\caption{10 connections with the highest heritability index for state 2. Connections are sorted with respect to HI values.}
\label{tab:my_table2}
\begin{tabular}{|l|l|c|}
\hline
Regions & Regions & HI \\
\hline
Hippocampus-R & Cerebelum-3-R & 1.00 \\
Hippocampus-L & Cerebelum-4-5-L & 0.93 \\
Occipital-Mid-R & Cerebelum-Crus2-R & 0.86 \\
Olfactory-L & Cerebelum-3-L & 0.81 \\
Heschl-L & Temporal-Pole-Sup-L & 0.81 \\
Rolandic-Oper-L & Temporal-Pole-Sup-L & 0.80 \\
Caudate-R & Cerebelum-Crus1-L & 0.79 \\
Cerebelum-7b-R & Cerebelum-9-R & 0.78 \\
Cingulum-Ant-L & Cerebelum-3-R & 0.78 \\
Frontal-Mid-Orb-R & Frontal-Med-Orb-L & 0.78 \\
\hline
\end{tabular}
\caption{10 connections with the highest heritability index for state 3. Connections are sorted with respect to HI values.}
\label{tab:my_table3}
\end{table}

\subsection{Null test on twin study design}
\begin{figure}[t]
\includegraphics[width=\linewidth]{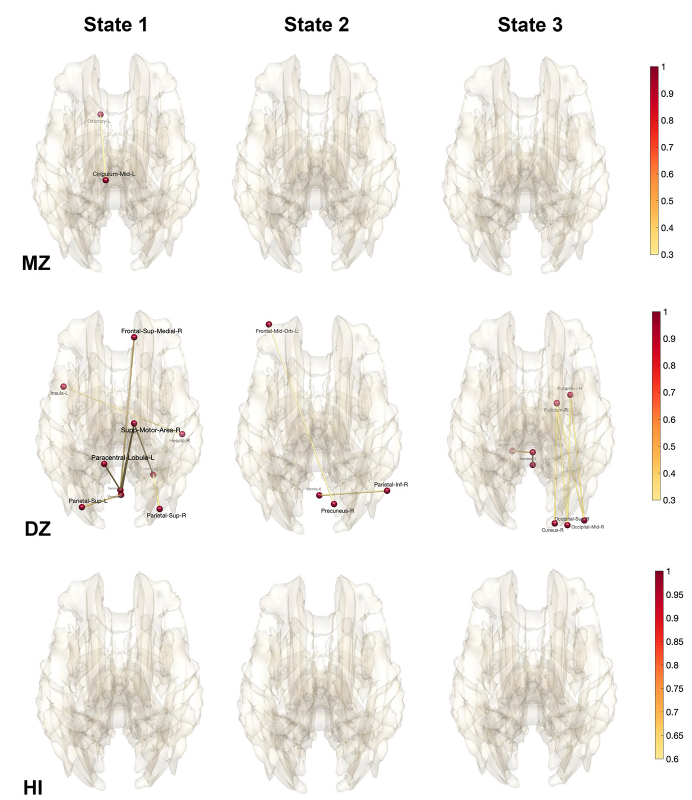}
\caption{\footnotesize MZ-correlation (top) and DZ-correlation (middle) in each state obtained through topological clustering in Figure \ref{fig:clustering-results}. There is no MZ-correlation above 0.3 and not displayed. The heritability index (HI) is determined by the twice the difference in twin correlations. HI of each state shows extensive genetic contribution of dynamically changing states. }
\label{fig:null}
\end{figure}

\begin{figure}[t]
\includegraphics[width=\linewidth]{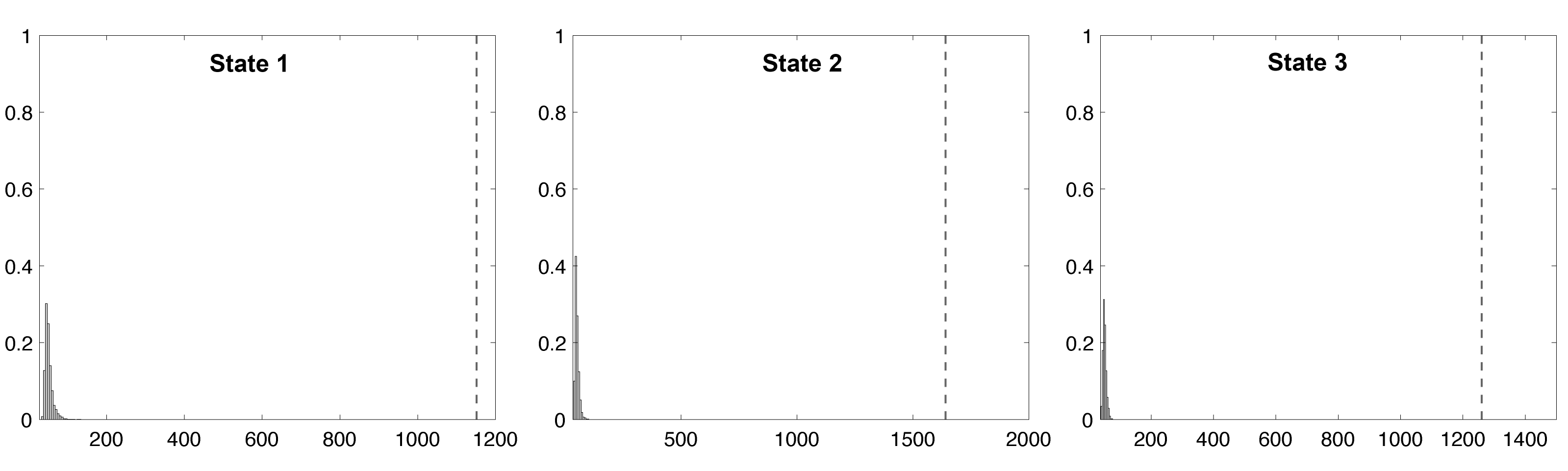}
\caption{\footnotesize The normalized histogram of  the Wasserstein distance between average MZ- and DZ-twin correlations within each state over 1000 derangements. Since the generated null data has no genetic signal, we are basically computing the Wasserstein distance between two connectivity matrices with random noises. In comparison,  the observed Wasserstein distance (dotted line) between average MZ- and DZ-twin correlation  shows huge topological differences.}
\label{fig:null2}
\end{figure}

Because we are reporting significantly higher diffused heritability compared to existing literature \citep{glahn.2010, xu.2017, korgaonkar.2014}, we performed the null test to check the validity of our analysis pipeline further. We generated the null MZ-twin data by randomly pairing each MZ individual with another, excluding their own twin. Such a permutation is generated by {\em derangement}, which is a permutation of the elements of a set, such that no element appears in its original position \citep{hassani.2003}. In other words, if we have a set of distinct items and you rearrange them, a derangement means none of the items are in the spot they started in. The null DZ-twin data is constructed similarly. Such null data should not show any genetic relations beyond random chances. On the null data, we recomputed the twin correlations and the heritability index,  following the same pipeline as before. Figure \ref{fig:null} shows an example of one possible derangement out of exponentially many such permutations. For $m$ MZ-twin pairs, there are $$m! \sum_{i=0}^m \frac{(-1)^i}{i!}$$ number of derangements. For the null test, we generated 1000 derangements then followed the proposed pipeline in computing average MZ- and DZ-correlations in each state.  We used the Wasserstein distance in measuring the topological discrepancy. Figure \ref{fig:null2} displays  the normalized histogram of  the Wasserstein distance between average MZ- and DZ-twin correlations within each state over 1000 derangements. Because the generated null data has no genetic signal, we are basically computing the Wasserstein distance between two random connectivity matrices. In comparison,  the observed Wasserstein distance (red line) between average MZ- and DZ-twin correlation  shows huge topological differences. None of the derangements show the large wide spread signals as our observation. We conclude that what we observe is  genetic signal and cannot possibly be produced by random chance.

\section{Discussion}

In this study, we proposed the topological clustering method for the estimation and quantification of dynamic state changes in time-varying brain networks. A coherent statistical theory, grounded in persistent homology, was developed, and we demonstrated the application of this method to resting-state fMRI data. Resting-state brain networks tend to persist in a single state for extended periods before transitioning to another state \citep{allen.2014,shakil.2016,calhoun.2016,rosenberg.2020}. The average brain network in each state  appears to diverge from the patterns reported in previous studies (Figure 10) \citep{cai.2018.TMI,haimovici.2017,huang.2020.NM}. Further research is required for independent validation.

In contrast to previous studies that reported relatively low heritability in functional brain networks \citep{glahn.2010, xu.2017, korgaonkar.2014, wan.2022}, our findings indicate significant higher heritability across various regions of the brain network. This discovery not only challenges the prevailing understanding but also opens new avenues for exploring genetic influences on brain network dynamics. Our observations align with the early findings by \citet{lykken.1982}, which documented higher heritability in EEG spectra. Heritability in brain networks may be more nuanced than previously understood. In our framework, rather than directly using the connectivity strength, we decomposed networks into discrete topological states and computed heritability for each state. This granular analysis provides a more accurate estimation of heritability across different functional states of the brain. The resting state measures employed in studies such as those by \citet{glahn.2010, xu.2017, korgaonkar.2014} directly rely on {\em static} connectivity matrices. These matrices, while informative, often do not capture the dynamic and configural nature of brain networks. Such methods may overlook hidden configural patterns that hold significant heritable information. Our topological method represents a significant advancement in this regard. By focusing on the topological aspects of dynamic brain networks, our method is adept at identifying and extracting hidden patterns of high heritability that might be missed by traditional approaches.  This capability could be crucial for understanding the genetic basis of various neuropsychiatric and neurodevelopmental disorders, where altered brain network configurations play a critical role.

Intraclass correlation (ICC) has long been recognized as a vital reliability and reproducibility metric, especially for gauging similarity in paired data when the order of pairing is not preserved \citep{chen.2018, sarker.2023, solis.2023}. In brain imaging, it serves as a popular baseline for test-retest (TRT) reliability assessments, often in conjunction with the Dice coefficient \citep{liao.2013, cole.2014, cousineau.2017, pfaehler.2021,zhang.2018}. The widespread use of ICC in these contexts underscores its perceived utility in evaluating consistency across imaging sessions or different imaging modalities. The conventional computation of ICC is typically through an ANOVA statistical model, which can be fairly limited and inflexible. Recent years have seen a shift towards mixed-effects models, which offer greater flexibility and accuracy in estimating ICC, especially in datasets with nested or hierarchical structures \citep{chen.2018,nielsen.2021}. In light of these advancements, our proposed transposition-based approach for computing correlation over paired data presents a novel approach to computing ICC, potentially offering a faster and more efficient alternative. The full potential and utility of the transposition-based method for ICC computation, however, remain to be explored in future research. 

\section*{Acknowledgement}
We would like to thank  Chee-Ming Ting and Hernando Ombao of  KAUST for discussion on $k$-means clustering. We also like to thank Soumya Das, Tananun Songdechakraiwut of University of Wisconsin, Madison and Botao Wang of University of Illinois, Urbana-Champaign for discussion on the clustering. 

\section*{Author Contributions}

{\bf Moo K. Chung}: Conceptualization, Formal Analysis, Funding Acquisition, Investigation, Methodology, Project Administration, Software, Supervision, Validation, Visualization, Writing- Original Draft Preparation, Writing- Review \& Editing\\

\noindent {\bf Shih-Gu Huang}: Data Curation, Formal Analysis, Investigation, Methodology, Software, Writing- Original Draft Preparation,  Writing- Review \& Editing\\

\noindent {\bf Ian C. Carrol}: Data Curation, Investigation, Resources, Software\\

\noindent {\bf Vince D. Calhoun}: Writing- Review \& Editing\\

\noindent {\bf  H. Hill Goldsmith}: Funding Acquisition, Project Administration, Resources, Supervision, Writing- Review \& Editing\\

\subsection*{Ethics statement}
The ethics approval for using the data  was obtained from the local Institutional Review Boards (IRB)  of University of Wisconsin-Madison (\url{https://irb.wisc.edu}).  Informed written consent was obtained from all participants.

 \subsection*{Data availability statement}
 The data used in the study is publicly available through National Data Archive (\url{https://nda.nih.gov}) under the collection tile Validating RDoC for Children and Adolescents: A Twin Study with Neuroimaging \#2105.
The data can be also obtained by contacting  Wisconsin Twin Research (\url{https://goldsmithtwins.waisman.wisc.edu}) through email {\tt wisconsintwins@waisman.wisc.edu}.

\section*{Funding}
This study was supported by the National Institutes of Health  (EB022856, MH133614 to MC; MH101504, P30HD003352, U54HD09025 to HG) and the National Science Foundation (MDS-2010778 to MKC; 2112455 to VC). The funders had no role in study design, data collection and analysis, decision to publish, or preparation of the manuscript.

\section*{Competing interests}
The authors have declared that no competing interests exist.

\bibliographystyle{plainnat}

\begin{thebibliography}{114}
\providecommand{\natexlab}[1]{#1}
\providecommand{\url}[1]{\texttt{#1}}
\expandafter\ifx\csname urlstyle\endcsname\relax
  \providecommand{\doi}[1]{doi: #1}\else
  \providecommand{\doi}{doi: \begingroup \urlstyle{rm}\Url}\fi

\bibitem[Agueh and Carlier(2011)]{agueh.2011}
M.~Agueh and G.~Carlier.
\newblock Barycenters in the Wasserstein space.
\newblock \emph{SIAM Journal on Mathematical Analysis}, 43:\penalty0 904--924,
  2011.

\bibitem[Aktas et~al.(2019)Aktas, Akbas, and Fatmaoui]{aktas.2019}
M.E. Aktas, E.~Akbas, and A.E. Fatmaoui.
\newblock Persistence homology of networks: methods and applications.
\newblock \emph{Applied Network Science}, 4:\penalty0 1--28, 2019.

\bibitem[Allen et~al.(2014)Allen, Damaraju, Plis, Erhardt, Eichele, and
  Calhoun]{allen.2014}
E.A. Allen, E.~Damaraju, S.M. Plis, E.B. Erhardt, T.~Eichele, and V.D. Calhoun.
\newblock Tracking whole-brain connectivity dynamics in the resting state.
\newblock \emph{Cerebral cortex}, 24:\penalty0 663--676, 2014.

\bibitem[Arbet et~al.(2020)Arbet, McGue, and Basu]{arbet.2020}
J.~Arbet, M.~McGue, and S.~Basu.
\newblock A robust and unified framework for estimating heritability in twin
  studies using generalized estimating equations.
\newblock \emph{Statistics in Medicine}, 2020.

\bibitem[Bassett and Sporns(2017)]{bassett.2017}
D.S. Bassett and O.~Sporns.
\newblock Network neuroscience.
\newblock \emph{Nature neuroscience}, 20\penalty0 (3):\penalty0 353--364, 2017.

\bibitem[Berwald et~al.(2018)Berwald, Gottlieb, and Munch]{berwald.2018}
J.J. Berwald, J.M. Gottlieb, and E.~Munch.
\newblock Computing wasserstein distance for persistence diagrams on a quantum
  computer.
\newblock \emph{arXiv:1809.06433}, 2018.

\bibitem[Billings et~al.(2021)Billings, Saggar, Hlinka, Keilholz, and
  Petri]{billings.2021}
J.~Billings, M.~Saggar, J.~Hlinka, S.~Keilholz, and G.~Petri.
\newblock Simplicial and topological descriptions of human brain dynamics.
\newblock \emph{Network Neuroscience}, 5:\penalty0 549--568, 2021.

\bibitem[Bishop(2006)]{bishop.2006}
C.M. Bishop.
\newblock \emph{Pattern recognition and machine learning}.
\newblock Springer, 2006.

\bibitem[Blokland et~al.(2011)Blokland, McMahon, Thompson, Martin,
  de~Zubicaray, and Wright]{blokland.2011}
G.A.M. Blokland, K.L. McMahon, P.M. Thompson, N.G. Martin, G.I. de~Zubicaray,
  and M.J. Wright.
\newblock Heritability of working memory brain activation.
\newblock \emph{The Journal of Neuroscience}, 31:\penalty0 10882--10890, 2011.

\bibitem[Bullmore and Sporns(2009)]{bullmore.2009}
E.~Bullmore and O.~Sporns.
\newblock {Complex brain networks: graph theoretical analysis of structural and
  functional systems}.
\newblock \emph{Nature Review Neuroscience}, 10:\penalty0 186--98, 2009.

\bibitem[Burghy et~al.(2016)Burghy, Fox, Cornejo, Stodola, Sommerfeldt,
  Westbrook, Van~Hulle, Schmidt, Goldsmith, Davidson, et~al.]{burghy.2016}
C.A. Burghy, M.E. Fox, M.D. Cornejo, D.E. Stodola, S.L. Sommerfeldt, C.A.
  Westbrook, C.~Van~Hulle, N.L. Schmidt, H.H. Goldsmith, R.J. Davidson, et~al.
\newblock Experience-driven differences in childhood cortisol predict
  affect-relevant brain function and coping in adolescent {M}onozygotic twins.
\newblock \emph{Scientific Reports}, 6:\penalty0 37081, 2016.

\bibitem[Cai et~al.(2018)Cai, Zille, Stephen, Wilson, Calhoun, and
  Wang]{cai.2018.TMI}
B.~Cai, P.~Zille, J.M. Stephen, T.W. Wilson, V.D. Calhoun, and Y.P. Wang.
\newblock Estimation of dynamic sparse connectivity patterns from resting state
  f{MRI}.
\newblock \emph{IEEE Transactions on Medical Imaging}, 37:\penalty0 1224--1234,
  2018.

\bibitem[Calhoun and Adali(2016)]{calhoun.2016}
V.D. Calhoun and T.~Adali.
\newblock Time-varying brain connectivity in f{MRI} data: whole-brain
  data-driven approaches for capturing and characterizing dynamic states.
\newblock \emph{IEEE Signal Processing Magazine}, 33:\penalty0 52--66, 2016.

\bibitem[Canas and Rosasco(2012)]{canas.2012}
G.D. Canas and L.~Rosasco.
\newblock Learning probability measures with respect to optimal transport
  metrics.
\newblock \emph{arXiv preprint arXiv:1209.1077}, 2012.

\bibitem[Caputi et~al.(2021)Caputi, Pidnebesna, and Hlinka]{caputi.2021}
L.~Caputi, A.~Pidnebesna, and J.~Hlinka.
\newblock Promises and pitfalls of topological data analysis for brain
  connectivity analysis.
\newblock \emph{NeuroImage}, 238:\penalty0 118245, 2021.

\bibitem[Chen et~al.(2019)Chen, Ni, Bai, and Wang]{chen.2019}
C.~Chen, X.~Ni, Q.~Bai, and Y.~Wang.
\newblock A topological regularizer for classifiers via persistent homology.
\newblock In \emph{The 22nd International Conference on Artificial Intelligence
  and Statistics}, pages 2573--2582. PMLR, 2019.

\bibitem[Chen et~al.(2018)Chen, Taylor, Haller, Kircanski, Stoddard, Pine,
  Leibenluft, Brotman, and Cox]{chen.2018}
G.~Chen, P.A. Taylor, S.P. Haller, K.~Kircanski, J.~Stoddard, D.S. Pine,
  E.~Leibenluft, M.A. Brotman, and R.W. Cox.
\newblock Intraclass correlation: Improved modeling approaches and applications
  for neuroimaging.
\newblock \emph{Human brain mapping}, 39:\penalty0 1187--1206, 2018.

\bibitem[Chiang et~al.(2011)Chiang, McMahon, de~Zubicaray, Martin, Hickie,
  Toga, Wright, and Thompson]{chiang.2011}
M.-C. Chiang, K.L. McMahon, G.I. de~Zubicaray, N.G. Martin, I.~Hickie, A.W.
  Toga, M.J. Wright, and P.M. Thompson.
\newblock Genetics of white matter development: a {DTI} study of 705 twins and
  their siblings aged 12 to 29.
\newblock \emph{NeuroImage}, 54:\penalty0 2308--2317, 2011.

\bibitem[Chung(2012)]{chung.2012.CNA}
M.K. Chung.
\newblock \emph{Computational Neuroanatomy: The Methods}.
\newblock World Scientific, Singapore, 2012.

\bibitem[Chung et~al.(2007)Chung, Dalton, Shen, Evans, and
  Davidson]{chung.2007.TMI}
M.K. Chung, K.M. Dalton, L.~Shen, A.C. Evans, and R.J. Davidson.
\newblock Weighted {Fourier} representation and its application to quantifying
  the amount of gray matter.
\newblock \emph{IEEE Transactions on Medical Imaging}, 26:\penalty0 566--581,
  2007.

\bibitem[Chung et~al.(2008)Chung, Dalton, and Davidson]{chung.2008.TMI}
M.K. Chung, K.M. Dalton, and R.J. Davidson.
\newblock Tensor-based cortical surface morphometry via weighted spherical
  harmonic representation.
\newblock \emph{IEEE Transactions on Medical Imaging}, 27:\penalty0 1143--1151,
  2008.

\bibitem[Chung et~al.(2013)Chung, Hanson, Lee, Adluru, Alexander, Davidson, and
  Pollak]{chung.2013.MICCAI}
M.K. Chung, J.L. Hanson, H.~Lee, Nagesh Adluru, Andrew~L. Alexander, R.J.
  Davidson, and S.D. Pollak.
\newblock Persistent homological sparse network approach to detecting white
  matter abnormality in maltreated children: {MRI} and {DTI} multimodal study.
\newblock \emph{MICCAI, Lecture Notes in Computer Science (LNCS)},
  8149:\penalty0 300--307, 2013.

\bibitem[Chung et~al.(2017{\natexlab{a}})Chung, Hanson, Adluru, Alexander,
  Davidson, and Pollak]{chung.2017.BC}
M.K. Chung, J.L. Hanson, L.~Adluru, A.L. Alexander, R.J. Davidson, and S.D.
  Pollak.
\newblock Integrative structural brain network analysis in diffusion tensor
  imaging.
\newblock \emph{Brain Connectivity}, 7:\penalty0 331--346, 2017{\natexlab{a}}.

\bibitem[Chung et~al.(2017{\natexlab{b}})Chung, Lee, Solo, Davidson, and
  Pollak]{chung.2017.CNI}
M.K. Chung, H.~Lee, V.~Solo, R.J. Davidson, and S.D. Pollak.
\newblock Topological distances between brain networks.
\newblock \emph{International Workshop on Connectomics in Neuroimaging},
  10511:\penalty0 161--170, 2017{\natexlab{b}}.

\bibitem[Chung et~al.(2019{\natexlab{a}})Chung, Huang, Gritsenko, Shen, and
  Lee]{chung.2019.ISBI}
M.K. Chung, S.-G. Huang, A.~Gritsenko, L.~Shen, and H.~Lee.
\newblock Statistical inference on the number of cycles in brain networks.
\newblock In \emph{2019 IEEE 16th International Symposium on Biomedical Imaging
  (ISBI 2019)}, pages 113--116. IEEE, 2019{\natexlab{a}}.

\bibitem[Chung et~al.(2019{\natexlab{b}})Chung, Lee, DiChristofano, Ombao, and
  Solo]{chung.2019.NN}
M.K. Chung, H.~Lee, A.~DiChristofano, H.~Ombao, and V.~Solo.
\newblock Exact topological inference of the resting-state brain networks in
  twins.
\newblock \emph{Network Neuroscience}, 3:\penalty0 674--694,
  2019{\natexlab{b}}.

\bibitem[Chung et~al.(2019{\natexlab{c}})Chung, Xie, Huang, Wang, Yan, and
  Shen]{chung.2019.CNI}
M.K. Chung, L.~Xie, S.-G. Huang, Y.~Wang, J.~Yan, and L.~Shen.
\newblock Rapid acceleration of the permutation test via transpositions.
\newblock 11848:\penalty0 42--53, 2019{\natexlab{c}}.

\bibitem[Chung et~al.(2023{\natexlab{a}})Chung, Das, and Ombao]{chung.2023.DE}
M.K. Chung, S.~Das, and H.~Ombao.
\newblock Topological data analysis of functional human brain networks.
\newblock \emph{arXiv preprint arXiv:2210.09092}, 2023{\natexlab{a}}.

\bibitem[Chung et~al.(2023{\natexlab{b}})Chung, Huang, Carroll, Calhoun, and
  Hill]{chung.2023.wasserstein}
M.K. Chung, S.-G. Huang, I.C. Carroll, V.D. Calhoun, and G.H Hill.
\newblock Persistent homological state-space estimation of functional human
  brain networks at rest.
\newblock \emph{arXiv e-prints}, pages arXiv--2201.00087, 2023{\natexlab{b}}.
\newblock URL \url{https://arxiv.org/pdf/2201.00087}.

\bibitem[Chung et~al.(2023{\natexlab{c}})Chung, Ramos, De~Paiva, Mathis,
  Prabhakaran, Nair, Meyerand, Hermann, Binder, and Struck]{chung.2023.NI}
M.K. Chung, C.G. Ramos, F.B. De~Paiva, J.~Mathis, V.~Prabhakaran, V.A. Nair,
  M.E. Meyerand, B.P. Hermann, J.R. Binder, and A.F. Struck.
\newblock Unified topological inference for brain networks in temporal lobe
  epilepsy using the {Wasserstein} distance.
\newblock \emph{NeuroImage}, 284:\penalty0 120436, 2023{\natexlab{c}}.

\bibitem[Cole et~al.(2014)Cole, Farmer, Rees, Johnson, Frost, Scahill, and
  Hobbs]{cole.2014}
J.H. Cole, R.E. Farmer, E.M. Rees, H.J. Johnson, C.~Frost, R.I. Scahill, and
  N.Z. Hobbs.
\newblock Test-retest reliability of diffusion tensor imaging in
  {HuntingtonÕs} disease.
\newblock \emph{PLoS Currents}, 6, 2014.

\bibitem[Cousineau et~al.(2017)Cousineau, Jodoin, Garyfallidis, C{\^o}t{\'e},
  Morency, Rozanski, GrandÕMaison, Bedell, and Descoteaux]{cousineau.2017}
M.~Cousineau, P.-M. Jodoin, E.~Garyfallidis, M.-A. C{\^o}t{\'e}, F.C. Morency,
  V.~Rozanski, M.~GrandÕMaison, B.J. Bedell, and M.~Descoteaux.
\newblock A test-retest study on {Parkinson's PPMI} dataset yields
  statistically significant white matter fascicles.
\newblock \emph{NeuroImage: Clinical}, 16:\penalty0 222--233, 2017.

\bibitem[Cuturi and Doucet(2014)]{cuturi.2014}
M.~Cuturi and A.~Doucet.
\newblock Fast computation of {W}asserstein barycenters.
\newblock In \emph{International conference on machine learning}, pages
  685--693. PMLR, 2014.

\bibitem[Devlin et~al.(1975)Devlin, Gnanadesikan, and Kettenring]{devlin.1975}
S.J. Devlin, R.~Gnanadesikan, and J.R. Kettenring.
\newblock Robust estimation and outlier detection with correlation
  coefficients.
\newblock \emph{Biometrika}, 62:\penalty0 531--545, 1975.

\bibitem[Dubey and M{\"u}ller(2019)]{dubey.2019}
P.~Dubey and H.-G. M{\"u}ller.
\newblock Fr{\'e}chet analysis of variance for random objects.
\newblock \emph{Biometrika}, 106:\penalty0 803--821, 2019.

\bibitem[Edelsbrunner and Harer(2010)]{edelsbrunner.2010}
H.~Edelsbrunner and J.~Harer.
\newblock \emph{Computational topology: {A}n introduction}.
\newblock American Mathematical Society, 2010.

\bibitem[Edmonds and Karp(1972)]{edmonds.1972}
J.~Edmonds and R.M. Karp.
\newblock Theoretical improvements in algorithmic efficiency for network flow
  problems.
\newblock \emph{Journal of the ACM (JACM)}, 19:\penalty0 248--264, 1972.

\bibitem[Falconer and Mackay(1995)]{falconer.1995}
D.~Falconer and T~Mackay.
\newblock \emph{Introduction to Quantitative Genetics, 4th ed.}
\newblock Longman, 1995.

\bibitem[Fu et~al.(2023)Fu, Huang, Zhang, Dong, Xue, Niu, Li, Shi, Wang, and
  Zhang]{fu.2023}
Y.~Fu, Y.~Huang, Z.~Zhang, S.~Dong, L.~Xue, M.~Niu, Y.~Li, Z.~Shi, Y.~Wang, and
  H.~Zhang.
\newblock {OTFPF}: Optimal transport based feature pyramid fusion network for
  brain age estimation.
\newblock \emph{Information Fusion}, 100:\penalty0 101931, 2023.

\bibitem[Ghrist(2008)]{ghrist.2008}
R.~Ghrist.
\newblock Barcodes: The persistent topology of data.
\newblock \emph{Bulletin of the American Mathematical Society}, 45:\penalty0
  61--75, 2008.

\bibitem[Giusti et~al.(2016)Giusti, Ghrist, and Bassett]{giusti.2016}
C.~Giusti, R.~Ghrist, and D.S. Bassett.
\newblock Two’s company, three (or more) is a simplex.
\newblock \emph{Journal of computational neuroscience}, 41:\penalty0 1--14,
  2016.

\bibitem[Glahn et~al.(2010)Glahn, Winkler, Kochunov, Almasy, Duggirala,
  Carless, Curran, Olvera, Laird, and Smith]{glahn.2010}
D.C. Glahn, A.M. Winkler, P.~Kochunov, L.~Almasy, R.~Duggirala, M.A. Carless,
  J.C. Curran, R.L. Olvera, A.R. Laird, and S.M. Smith.
\newblock Genetic control over the resting brain.
\newblock \emph{Proceedings of the National Academy of Sciences}, 107:\penalty0
  1223--1228, 2010.

\bibitem[Gupta et~al.(2022)Gupta, Hu, Kaan, Jin, Mpoy, Chung, Singh, Saltz,
  Kurc, Saltz, et~al.]{gupta.2022}
S.~Gupta, X.~Hu, J.~Kaan, M.~Jin, M.~Mpoy, K.~Chung, G.~Singh, M.~Saltz,
  T.~Kurc, J.~Saltz, et~al.
\newblock Learning topological interactions for multi-class medical image
  segmentation.
\newblock In \emph{European Conference on Computer Vision}, pages 701--718,
  2022.

\bibitem[Haimovici et~al.(2017)Haimovici, Tagliazucchi, Balenzuela, and
  Laufs]{haimovici.2017}
A.~Haimovici, E.~Tagliazucchi, P.~Balenzuela, and H.~Laufs.
\newblock On wakefulness fluctuations as a source of {BOLD} functional
  connectivity dynamics.
\newblock \emph{Scientific Reports}, 7:\penalty0 5908, 2017.

\bibitem[Hartigan and Wong(1979)]{hartigan.1979}
J.A. Hartigan and M.A. Wong.
\newblock Algorithm {AS} 136: A k-means clustering algorithm.
\newblock \emph{Journal of the Royal Statistical Society. Series C (applied
  statistics)}, 28:\penalty0 100--108, 1979.

\bibitem[Hartmann et~al.(2018)Hartmann, Schirrmeister, and Ball]{hartmann.2018}
K.G. Hartmann, R.T. Schirrmeister, and T.~Ball.
\newblock EEG-GAN: Generative adversarial networks for electroencephalograhic
  (EEG) brain signals.
\newblock \emph{arXiv preprint arXiv:1806.01875}, 2018.

\bibitem[Hassani(2003)]{hassani.2003}
M.~Hassani.
\newblock Derangements and applications.
\newblock \emph{Journal of Integer Sequences}, 6:\penalty0 03--1, 2003.

\bibitem[Hofer et~al.(2019)Hofer, Kwitt, Niethammer, and Dixit]{hofer.2019}
C.~Hofer, R.~Kwitt, M.~Niethammer, and M.~Dixit.
\newblock Connectivity-optimized representation learning via persistent
  homology.
\newblock In \emph{International Conference on Machine Learning}, pages
  2751--2760, 2019.

\bibitem[Hu et~al.(2019)Hu, Li, Samaras, and Chen]{hu.2019}
X.~Hu, F.~Li, D.~Samaras, and C.~Chen.
\newblock Topology-preserving deep image segmentation.
\newblock In \emph{Advances in Neural Information Processing Systems}, pages
  5657--5668, 2019.

\bibitem[Huang et~al.(2019{\natexlab{a}})Huang, Chung, Carroll, and
  Goldsmith]{huang.2019.DSW}
S.-G. Huang, M.~K. Chung, I.~C. Carroll, and H.~H. Goldsmith.
\newblock Dynamic functional connectivity using heat kernel.
\newblock In \emph{2019 IEEE Data Science Workshop (DSW)}, pages 222--226,
  2019{\natexlab{a}}.
\newblock \doi{10.1109/DSW.2019.8755550}.

\bibitem[Huang et~al.(2019{\natexlab{b}})Huang, Gritsenko, Lindquist, and
  Chung]{huang.2019.ISBI}
S.-G. Huang, A.~Gritsenko, M.A. Lindquist, and M.K. Chung.
\newblock Circular pearson correlation using cosine series expansion.
\newblock In \emph{IEEE 16th International Symposium on Biomedical Imaging
  (ISBI)}, pages 1774--1777, 2019{\natexlab{b}}.

\bibitem[Huang et~al.(2020)Huang, Samdin, Ting, Ombao, and
  Chung]{huang.2020.NM}
S.-G. Huang, S.-T. Samdin, C.M. Ting, H.~Ombao, and M.K. Chung.
\newblock Statistical model for dynamically-changing correlation matrices with
  application to brain connectivity.
\newblock \emph{Journal of Neuroscience Methods}, 331:\penalty0 108480, 2020.

\bibitem[Hutchison et~al.(2013)Hutchison, Womelsdorf, Allen, Bandettini, and
  Calhoun]{hutchison.2013}
R.M. Hutchison, T.~Womelsdorf, E.A. Allen, P.A. Bandettini, and V.D. et.~al.
  Calhoun.
\newblock Dynamic functional connectivity: promise, issues, and
  interpretations.
\newblock \emph{NeuroImage}, 80:\penalty0 360--378, 2013.

\bibitem[Jenkinson et~al.(2002)Jenkinson, Bannister, Brady, and
  Smith]{jenkinson.2002}
M.~Jenkinson, P.~Bannister, M.~Brady, and S.~Smith.
\newblock {Improved optimization for the robust and accurate linear
  registration and motion correction of brain images}.
\newblock \emph{NeuroImage}, 17:\penalty0 825--841, 2002.

\bibitem[Johnson(1967)]{johnson.1967}
S.C. Johnson.
\newblock Hierarchical clustering schemes.
\newblock \emph{Psychometrika}, 32:\penalty0 241--254, 1967.

\bibitem[Khalid et~al.(2014)Khalid, Kim, Chung, Ye, and Jeon]{khalid.2014}
A.~Khalid, B.S. Kim, M.K. Chung, J.C. Ye, and D.~Jeon.
\newblock Tracing the evolution of multi-scale functional networks in a mouse
  model of depression using persistent brain network homology.
\newblock \emph{NeuroImage}, 101:\penalty0 351--363, 2014.

\bibitem[Korgaonkar et~al.(2014)Korgaonkar, Ram, Williams, Gatt, and
  Grieve]{korgaonkar.2014}
M.S. Korgaonkar, K.~Ram, L.M. Williams, J.M. Gatt, and S.M. Grieve.
\newblock Establishing the resting state default mode network derived from
  functional magnetic resonance imaging tasks as an endophenotype: a twins
  study.
\newblock \emph{Human brain mapping}, 35:\penalty0 3893--3902, 2014.

\bibitem[Kuang et~al.(2019)Kuang, Zhao, Xing, Chen, Xiong, and Han]{kuang.2019}
Liqun Kuang, Deyu Zhao, Jiacheng Xing, Zhongyu Chen, Fengguang Xiong, and Xie
  Han.
\newblock Metabolic brain network analysis of fdg-pet in alzheimer’s disease
  using kernel-based persistent features.
\newblock \emph{Molecules}, 24\penalty0 (12):\penalty0 2301, 2019.

\bibitem[Kullback and Leibler(1951)]{kullback.1951}
S.~Kullback and R.A. Leibler.
\newblock On information and sufficiency.
\newblock \emph{The Annals of Mathematical Statistics}, 22:\penalty0 79--86,
  1951.

\bibitem[Le and Kume(2000)]{le.2000}
H.~Le and A.~Kume.
\newblock The {Fr{\'e}chet} mean shape and the shape of the means.
\newblock \emph{Advances in Applied Probability}, 32:\penalty0 101--113, 2000.

\bibitem[Lee et~al.(2011)Lee, Chung, Kang, Kim, and Lee]{lee.2011.MICCAI}
H.~Lee, M.K. Chung, H.~Kang, B.-N. Kim, and D.S. Lee.
\newblock Computing the shape of brain networks using graph filtration and
  {Gromov-Hausdorff} metric.
\newblock \emph{{MICCAI, Lecture Notes in Computer Science}}, 6892:\penalty0
  302--309, 2011.

\bibitem[Lee et~al.(2012)Lee, Kang, Chung, Kim, and Lee]{lee.2012.tmi}
H.~Lee, H.~Kang, M.K. Chung, B.-N. Kim, and D.S Lee.
\newblock Persistent brain network homology from the perspective of dendrogram.
\newblock \emph{IEEE Transactions on Medical Imaging}, 31:\penalty0 2267--2277,
  2012.

\bibitem[Liao et~al.(2013)Liao, Xia, Xu, Dai, Cao, Niu, Zuo, Zang, and
  He]{liao.2013}
X.-H. Liao, M.-R. Xia, T.~Xu, Z.-J. Dai, X.-Y. Cao, H.-J. Niu, X.-N. Zuo, Y.-F.
  Zang, and Y.~He.
\newblock Functional brain hubs and their test--retest reliability: a multiband
  resting-state functional MRI study.
\newblock \emph{NeuroImage}, 83:\penalty0 969--982, 2013.

\bibitem[Lin et~al.(2023)Lin, Zepf, Christensen, Bashir, Svendsen, Tolsgaard,
  and Feragen]{lin.2023}
M.~Lin, K.~Zepf, A.N. Christensen, Z.~Bashir, M.B.S. Svendsen, M.~Tolsgaard,
  and A.~Feragen.
\newblock {DTU-Net}: Learning topological similarity for curvilinear structure
  segmentation.
\newblock In \emph{International Conference on Information Processing in
  Medical Imaging}, pages 654--666, 2023.

\bibitem[Lindquist(2014)]{lindquist.2014}
M~Lindquist.
\newblock Statistical and computational methods in brain image analysis. by
  {Moo K. Chung}. {Boca Raton}, {Florida}: {CRC} press. 2013.
\newblock \emph{Journal of the American Statistical Association}, 109:\penalty0
  1334--1335, 2014.

\bibitem[Lykken et~al.(1982)Lykken, Tellegen, and Iacono]{lykken.1982}
D.T. Lykken, A.~Tellegen, and W.G. Iacono.
\newblock {EEG} spectra in twins: Evidence for a neglected mechanism of genetic
  determination.
\newblock \emph{Physiological Psychology}, 10:\penalty0 60--65, 1982.

\bibitem[Ma et~al.(2023)Ma, Wen, Zhu, and Zhang]{ma.2023}
K.~Ma, X.~Wen, Q.~Zhu, and D.~Zhang.
\newblock Positive definite wasserstein graph kernel for brain disease
  diagnosis.
\newblock In \emph{International Conference on Medical Image Computing and
  Computer-Assisted Intervention}, pages 168--177, 2023.

\bibitem[McKay et~al.(2014)McKay, Knowles, Winkler, Sprooten, Kochunov, Olvera,
  Curran, Kent~Jr., Carless, G{\"o}ring, Dyer, Duggirala, Almasy, Fox,
  Blangero, and Glahn]{mckay.2014}
D.R. McKay, E.E.M. Knowles, A.A.M. Winkler, E.~Sprooten, P.~Kochunov, R.L.
  Olvera, J.E. Curran, J.W. Kent~Jr., M.A. Carless, H.H.H. G{\"o}ring, T.D.
  Dyer, R.~Duggirala, L.~Almasy, P.T. Fox, J.~Blangero, and D.C. Glahn.
\newblock Influence of age, sex and genetic factors on the human brain.
\newblock \emph{Brain Imaging and Behavior}, 8:\penalty0 143--152, 2014.

\bibitem[Mi et~al.(2018)Mi, Zhang, Gu, and Wang]{mi.2018}
L.~Mi, W.~Zhang, X.~Gu, and Y.~Wang.
\newblock Variational wasserstein clustering.
\newblock In \emph{Proceedings of the European Conference on Computer Vision
  (ECCV)}, pages 322--337, 2018.

\bibitem[Mokhtari et~al.(2019)Mokhtari, Akhlaghi, Simpson, Wu, and
  Laurienti]{mokhtari.2019}
F.~Mokhtari, M.I. Akhlaghi, S.L. Simpson, G.~Wu, and P.J. Laurienti.
\newblock Sliding window correlation analysis: Modulating window shape for
  dynamic brain connectivity in resting state.
\newblock \emph{NeuroImage}, 189:\penalty0 655--666, 2019.

\bibitem[Nielsen et~al.(2021)Nielsen, Smink, and Fox]{nielsen.2021}
N.M. Nielsen, W.A.C. Smink, and J.-P. Fox.
\newblock Small and negative correlations among clustered observations:
  Limitations of the linear mixed effects model.
\newblock \emph{Behaviormetrika}, 48:\penalty0 51--77, 2021.

\bibitem[Oppenheim et~al.(1999)Oppenheim, Schafer, and Buck]{oppenheim.2001}
A.V. Oppenheim, R.W. Schafer, and J.R. Buck.
\newblock \emph{Discrete-time signal processing}.
\newblock Upper Saddle River, NJ: Prentice Hall, 1999.

\bibitem[Petri et~al.(2014)Petri, Expert, Turkheimer, Carhart-Harris, Nutt,
  Hellyer, and Vaccarino]{petri.2014}
G.~Petri, P.~Expert, F.~Turkheimer, R.~Carhart-Harris, D.~Nutt, P.J. Hellyer,
  and F.~Vaccarino.
\newblock Homological scaffolds of brain functional networks.
\newblock \emph{Journal of The Royal Society Interface}, 11:\penalty0 20140873,
  2014.

\bibitem[Pfaehler et~al.(2021)Pfaehler, Mesotten, Kramer, Thomeer, Vanhove,
  de~Jong, Adriaensens, Hoekstra, and Boellaard]{pfaehler.2021}
E.~Pfaehler, L.~Mesotten, G.~Kramer, M.~Thomeer, K.~Vanhove, J.~de~Jong,
  P.~Adriaensens, O.~S Hoekstra, and R.~Boellaard.
\newblock Repeatability of two semi-automatic artificial intelligence
  approaches for tumor segmentation in {PET}.
\newblock \emph{EJNMMI research}, 11:\penalty0 1--11, 2021.

\bibitem[Pozzi et~al.(2012)Pozzi, Di~Matteo, and Aste]{pozzi.2012}
F.~Pozzi, T.~Di~Matteo, and T.~Aste.
\newblock Exponential smoothing weighted correlations.
\newblock \emph{The European Physical Journal B}, 85:\penalty0 1--21, 2012.

\bibitem[Rabin et~al.(2011)Rabin, Peyr{\'e}, Delon, and Bernot]{rabin.2011}
J.~Rabin, G.~Peyr{\'e}, J.~Delon, and M.~Bernot.
\newblock Wasserstein barycenter and its application to texture mixing.
\newblock In \emph{International Conference on Scale Space and Variational
  Methods in Computer Vision}, pages 435--446. Springer, 2011.

\bibitem[Rashid et~al.(2014)Rashid, Damaraju, Pearlson, and
  Calhoun]{rashid.2014}
B.~Rashid, E.~Damaraju, G.D. Pearlson, and V.D. Calhoun.
\newblock Dynamic connectivity states estimated from resting f{MRI} identify
  differences among schizophrenia, bipolar disorder, and healthy control
  subjects.
\newblock \emph{Frontiers in Human Neuroscience}, 8:\penalty0 897, 2014.

\bibitem[Reynolds and Phillips(2015)]{reynolds.2015}
C.A. Reynolds and D.~Phillips.
\newblock Genetics of brain aging--twin aging.
\newblock 2015.

\bibitem[Rosenberg et~al.(2020)Rosenberg, Mennigen, Monti, and
  Kaiser]{rosenberg.2020}
B.M. Rosenberg, E.~Mennigen, M.M. Monti, and R.H. Kaiser.
\newblock Functional segregation of human brain networks across the lifespan:
  an exploratory analysis of static and dynamic resting-state functional
  connectivity.
\newblock \emph{Frontiers in Neuroscience}, 14:\penalty0 561594, 2020.

\bibitem[Sabbagh et~al.(2019)Sabbagh, Ablin, Varoquaux, Gramfort, and
  Engemann]{sabbagh.2019}
D.~Sabbagh, P.~Ablin, G.~Varoquaux, A.~Gramfort, and D.A. Engemann.
\newblock Manifold-regression to predict from meg/eeg brain signals without
  source modeling.
\newblock \emph{arXiv preprint arXiv:1906.02687}, 2019.

\bibitem[Sahu and Prasuna(2016)]{sahu.2016}
M.~Sahu and J.G. Prasuna.
\newblock Twin studies: A unique epidemiological tool.
\newblock \emph{Indian journal of community medicine: official publication of
  Indian Association of Preventive \& Social Medicine}, 41:\penalty0 177, 2016.

\bibitem[Santos et~al.(2019)Santos, Raposo, Coutinho-Filho, Copelli, Stam, and
  Douw]{santos.2019}
F.A.N. Santos, E.P. Raposo, M.D. Coutinho-Filho, M.~Copelli, C.J. Stam, and
  L.~Douw.
\newblock Topological phase transitions in functional brain networks.
\newblock \emph{Physical Review E}, 100:\penalty0 032414, 2019.

\bibitem[Sarker et~al.(2023)Sarker, Zaman, Ong, Paladugu, Aldred, Waisberg,
  Lee, and Tavakkoli]{sarker.2023}
P.~Sarker, N.~Zaman, J.~Ong, P.~Paladugu, M.~Aldred, E.~Waisberg, A.G. Lee, and
  A.~Tavakkoli.
\newblock Test--retest reliability of virtual reality devices in quantifying
  for relative afferent pupillary defect.
\newblock \emph{Translational Vision Science \& Technology}, 12:\penalty0 2,
  2023.

\bibitem[Shakil et~al.(2016)Shakil, Lee, and Keilholz]{shakil.2016}
S.~Shakil, C.-H. Lee, and S.D. Keilholz.
\newblock Evaluation of sliding window correlation performance for
  characterizing dynamic functional connectivity and brain states.
\newblock \emph{NeuroImage}, 133:\penalty0 111--128, 2016.

\bibitem[Shi et~al.(2016)Shi, Zhang, and Wang]{shi.2016}
J.~Shi, W.~Zhang, and Y.~Wang.
\newblock Shape analysis with hyperbolic wasserstein distance.
\newblock In \emph{Proceedings of the IEEE conference on computer vision and
  pattern recognition}, pages 5051--5061, 2016.

\bibitem[Sizemore et~al.(2018)Sizemore, Giusti, Kahn, Vettel, Betzel, and
  Bassett]{sizemore.2018}
A.E. Sizemore, C.~Giusti, A.~Kahn, J.M. Vettel, R.F. Betzel, and D.S. Bassett.
\newblock Cliques and cavities in the human connectome.
\newblock \emph{Journal of computational neuroscience}, 44:\penalty0 115--145,
  2018.

\bibitem[Sizemore et~al.(2019)Sizemore, Phillips-Cremins, Ghrist, and
  Bassett]{sizemore.2019}
A.E. Sizemore, J.E. Phillips-Cremins, R.~Ghrist, and D.S. Bassett.
\newblock The importance of the whole: topological data analysis for the
  network neuroscientist.
\newblock \emph{Network Neuroscience}, 3:\penalty0 656--673, 2019.

\bibitem[Smit et~al.(2008)Smit, Stam, Posthuma, Boomsma, and
  De~Geus]{smit.2008}
D.J.A. Smit, C.J. Stam, D.~Posthuma, D.I. Boomsma, and E.J.C. De~Geus.
\newblock Heritability of small-world networks in the brain: a graph
  theoretical analysis of resting-state {EEG} functional connectivity.
\newblock \emph{Human Brain Mapping}, 29:\penalty0 1368--1378, 2008.

\bibitem[Sol{\'\i}s-Lemus et~al.(2023)Sol{\'\i}s-Lemus, Baptiste, Barrows,
  Sillett, Gharaviri, Raffaele, Razeghi, Strocchi, Sim, and
  Kotadia]{solis.2023}
J.A. Sol{\'\i}s-Lemus, T.~Baptiste, R.~Barrows, C.~Sillett, A.~Gharaviri,
  G.~Raffaele, O.~Razeghi, M.~Strocchi, I.~Sim, and I.~Kotadia.
\newblock Evaluation of an open-source pipeline to create patient-specific left
  atrial models: A reproducibility study.
\newblock \emph{Computers in Biology and Medicine}, 162:\penalty0 107009, 2023.

\bibitem[Songdechakraiwut and Chung(2020)]{song.2020.ISBI}
T.~Songdechakraiwut and M.K Chung.
\newblock Dynamic topological data analysis for functional brain signals.
\newblock \emph{IEEE International Symposium on Biomedical Imaging Workshops},
  1:\penalty0 1--4, 2020.

\bibitem[Songdechakraiwut and Chung(2023)]{song.2023}
T.~Songdechakraiwut and M.K Chung.
\newblock Topological learning for brain networks.
\newblock \emph{Annals of Applied Statistics}, 17:\penalty0 403--433, 2023.

\bibitem[Songdechakraiwut et~al.(2021)Songdechakraiwut, Shen, and
  Chung]{song.2021.MICCAI}
T.~Songdechakraiwut, L.~Shen, and M.K. Chung.
\newblock Topological learning and its application to multimodal brain network
  integration.
\newblock \emph{Medical Image Computing and Computer Assisted Intervention
  (MICCAI)}, 12902:\penalty0 166--176, 2021.

\bibitem[Sporns(2003)]{sporns.2003}
O.~Sporns.
\newblock \emph{Graph Theory Methods for the Analysis of Neural Connectivity
  Patterns}, pages 171--185.
\newblock Springer US, Boston, MA, 2003.

\bibitem[Su et~al.(2015)Su, Zeng, Wang, Lu, and Gu]{su.2015}
Z.~Su, W.~Zeng, Y.~Wang, Z.-L. Lu, and X.~Gu.
\newblock Shape classification using wasserstein distance for brain morphometry
  analysis.
\newblock In \emph{International Conference on Information Processing in
  Medical Imaging}, pages 411--423. Springer, 2015.

\bibitem[Ting et~al.(2018)Ting, Ombao, Samdin, and Salleh]{ting.2018}
C.-M. Ting, H.~Ombao, S.B. Samdin, and S.-H. Salleh.
\newblock Estimating dynamic connectivity states in {fMRI} using
  regime-switching factor models.
\newblock \emph{IEEE transactions on Medical imaging}, 37:\penalty0 1011--1023,
  2018.

\bibitem[Turner et~al.(2014)Turner, Mileyko, Mukherjee, and Harer]{turner.2014}
K.~Turner, Y.~Mileyko, S.~Mukherjee, and J.~Harer.
\newblock Fr{\'e}chet means for distributions of persistence diagrams.
\newblock \emph{Discrete \& Computational Geometry}, 52:\penalty0 44--70, 2014.

\bibitem[Tzourio-Mazoyer et~al.(2002)Tzourio-Mazoyer, Landeau, Papathanassiou,
  Crivello, Etard, Delcroix, Mazoyer, and Joliot]{tzourio.2002}
N.~Tzourio-Mazoyer, B.~Landeau, D.~Papathanassiou, F.~Crivello, O.~Etard,
  N.~Delcroix, B.~Mazoyer, and M.~Joliot.
\newblock Automated anatomical labeling of activations in spm using a
  macroscopic anatomical parcellation of the {MNI MRI} single-subject brain.
\newblock \emph{NeuroImage}, 15:\penalty0 273--289, 2002.

\bibitem[Vaccarino et~al.(2022)Vaccarino, Fugacci, and
  Scaramuccia]{vaccarino.2022}
F.~Vaccarino, U.~Fugacci, and S.~Scaramuccia.
\newblock Persistent homology: A topological tool for higher-interaction
  systems.
\newblock In \emph{Higher-Order Systems}, pages 97--139. 2022.

\bibitem[Vallender(1974)]{vallender.1974}
S.S. Vallender.
\newblock Calculation of the {W}asserstein distance between probability
  distributions on the line.
\newblock \emph{Theory of Probability \& Its Applications}, 18:\penalty0
  784--786, 1974.

\bibitem[Vidaurre et~al.(2017)Vidaurre, Smith, and Woolrich]{vidaurre.2017}
D.~Vidaurre, S.M. Smith, and M.W. Woolrich.
\newblock Brain network dynamics are hierarchically organized in time.
\newblock \emph{Proceedings of the National Academy of Sciences}, 114:\penalty0
  12827--12832, 2017.

\bibitem[Wan et~al.(2022)Wan, Bayrak, Xu, Schaare, Bethlehem, Bernhardt, and
  Valk]{wan.2022}
B.~Wan, {\c{S}}.~Bayrak, T.~Xu, H.L. Schaare, R.~Bethlehem, B.C. Bernhardt, and
  S.L. Valk.
\newblock Heritability and cross-species comparisons of human cortical
  functional organization asymmetry.
\newblock \emph{Elife}, 11:\penalty0 e77215, 2022.

\bibitem[Wang et~al.(2017)Wang, Chung, Dentico, Lutz, and
  Davidson]{wang.2017.CNI}
Y.~Wang, M.K. Chung, D.~Dentico, A.~Lutz, and R.J. Davidson.
\newblock Topological network analysis of electroencephalographic power maps.
\newblock In \emph{International Workshop on Connectomics in NeuroImaging,
  Lecture Notes in Computer Science (LNCS)}, volume 10511, pages 134--142,
  2017.

\bibitem[Wang et~al.(2018)Wang, Ombao, and Chung]{wang.2018.annals}
Y.~Wang, H.~Ombao, and M.K. Chung.
\newblock Topological data analysis of single-trial electroencephalographic
  signals.
\newblock \emph{Annals of Applied Statistics}, 12:\penalty0 1506--1534, 2018.

\bibitem[Wijk et~al.(2010)Wijk, Stam, and Daffertshofer]{vanwijk.2010}
B.~C.~M. Wijk, C.~J. Stam, and A.~Daffertshofer.
\newblock Comparing brain networks of different size and connectivity density
  using graph theory.
\newblock \emph{PloS one}, 5:\penalty0 e13701, 2010.

\bibitem[Xing et~al.(2022)Xing, Jia, Wu, and Kuang]{xing.2022}
J.~Xing, J.~Jia, X.~Wu, and L.~Kuang.
\newblock A spatiotemporal brain network analysis of Alzheimer’s disease
  based on persistent homology.
\newblock \emph{Frontiers in aging neuroscience}, 14:\penalty0 788571, 2022.

\bibitem[Xu et~al.(2017)Xu, Yin, Ge, Han, Pang, Liu, Liu, and Friston]{xu.2017}
J.~Xu, X.~Yin, H.~Ge, Y.~Han, Z.~Pang, B.~Liu, S.~Liu, and K.~Friston.
\newblock Heritability of the effective connectivity in the resting-state
  default mode network.
\newblock \emph{Cerebral Cortex}, 27:\penalty0 5626--5634, 2017.

\bibitem[Xu et~al.(2021)Xu, Sanz, Garces, Maestu, Li, and Pantazis]{xu.2021}
Mengjia Xu, David~Lopez Sanz, Pilar Garces, Fernando Maestu, Quanzheng Li, and
  Dimitrios Pantazis.
\newblock A graph {G}aussian embedding method for predicting {A}lzheimer's
  disease progression with {MEG} brain networks.
\newblock \emph{IEEE Transactions on Biomedical Engineering}, 68:\penalty0
  1579--1588, 2021.

\bibitem[Yang et~al.(2020)Yang, Wen, and Davatzikos]{yang.2020}
Z.~Yang, J.~Wen, and C.~Davatzikos.
\newblock Smile-{GANs}: Semi-supervised clustering via {GANs} for dissecting
  brain disease heterogeneity from medical images.
\newblock \emph{arXiv preprint}, arXiv:\penalty0 2006.15255, 2020.

\bibitem[Yoo et~al.(2016)Yoo, Kim, Ahn, and Ye]{yoo.2016}
J.~Yoo, E.Y. Kim, Y.M. Ahn, and J.C. Ye.
\newblock Topological persistence vineyard for dynamic functional brain
  connectivity during resting and gaming stages.
\newblock \emph{Journal of neuroscience methods}, 267:\penalty0 1--13, 2016.

\bibitem[Yoo et~al.(2017)Yoo, Lee, Chung, Sohn, Chung, Na, Ju, and
  Jeong]{yoo.2017}
K.~Yoo, P.~Lee, M.K. Chung, W.S. Sohn, S.J. Chung, D.L. Na, D.~Ju, and
  Y.~Jeong.
\newblock Degree-based statistic and center persistency for brain connectivity
  analysis.
\newblock \emph{Human Brain Mapping}, 38:\penalty0 165--181, 2017.

\bibitem[Zemel and Panaretos(2019)]{zemel.2019}
Y.~Zemel and V.M. Panaretos.
\newblock Fr{\'e}chet means and procrustes analysis in {Wasserstein} space.
\newblock \emph{Bernoulli}, 25:\penalty0 932--976, 2019.

\bibitem[Zhan et~al.(2022)Zhan, Nagesh, Dean, Alexander, and
  Goldsmith]{zhan.2022}
L.~Zhan, A.~Nagesh, D.C. Dean, A.L. Alexander, and H.H. Goldsmith.
\newblock Genetic and environmental influences of variation in diffusion {MRI}
  measures of white matter microstructure.
\newblock \emph{Brain Structure and Function}, 227:\penalty0 131--144, 2022.

\bibitem[Zhang et~al.(2018)Zhang, Descoteaux, Zhang, Girard, Chamberland,
  Dunson, Srivastava, and Zhu]{zhang.2018}
Z.~Zhang, M.~Descoteaux, J.~Zhang, G.~Girard, M.~Chamberland, D.~Dunson,
  A.~Srivastava, and H.~Zhu.
\newblock Mapping population-based structural connectomes.
\newblock \emph{NeuroImage}, 172:\penalty0 130--145, 2018.

\bibitem[Zomorodian(2009)]{zomorodian.2009}
A.J. Zomorodian.
\newblock \emph{Topology for computing}.
\newblock Cambridge University Press, Cambridge, 2009.

\end{thebibliography}

\end{document}